\documentclass[english]{article}
\usepackage[T1]{fontenc}
\usepackage[latin9]{inputenc}
\usepackage[letterpaper]{geometry}
\usepackage{float}
\usepackage{comment}
\usepackage{amsmath}
\usepackage{graphicx}
\usepackage{color}
\usepackage{epstopdf}
\usepackage{adjustbox}
\usepackage{multirow}
\usepackage{array}
\usepackage{mathrsfs}
\usepackage{amsmath, amssymb, amsthm}
\usepackage{subfigure}
\usepackage{tikz}

\graphicspath{{./Figures/}}

\newtheorem{theorem}{Theorem}

\newtheorem{lemma}{Lemma}
\newtheorem{remark}{Remark}

\usepackage{amsfonts} 
\usepackage{color}
\definecolor{black}{rgb}{0,0,0}

\definecolor{red}{rgb}{1,0,0}

\definecolor{blue}{rgb}{0,0,1}

\newcommand{\specialcell}[2][c]{\begin{tabular}[#1]{@{}c@{}}#2\end{tabular}}
\usepackage{multirow}

\makeatother

\usepackage{babel}

\title{Residual driven  online mortar mixed finite element methods
and applications}

\author {Yanfang Yang\thanks{School of Mathematics and Information Science, Guangzhou University, Guangzhou, People's Republic of China.} \and
Eric T. Chung\thanks{Department of Mathematics, The Chinese University of Hong Kong, Hong Kong SAR. } \and
	Shubin Fu\thanks{Department of Mathematics, Texas A\&M University, College Station, TX 77843.}}

\begin{document}
		\maketitle
\begin{abstract}
	In this paper, we develop an online  basis enrichment method with the mortar mixed finite element method, using the oversampling technique, to solve for 
	flow problems in highly heterogeneous media. We first compute a coarse grid
	 solution with  a certain number of offline basis functions per edge, which are chosen as standard polynomials basis functions. We then iteratively enrich the multiscale solution space with online multiscale basis functions computed by using residuals. The iterative solution  converges to
	the fine scale solution rapidly. We also propose  an oversampling online method to achieve
	faster convergence speed.  The oversampling refers to using larger local regions in computing the online multiscale basis functions. We present extensive numerical experiments(including both 2D and 3D) to demonstrate the
	performance of our methods for both steady state flow, and two-phase flow and transport problems. In particular, for the time dependent two-phase flow and transport problems,  we apply the online method to the initial model,  without updating basis along the time evolution. Our numerical results demonstrate that by using a few number of  online basis functions, one can achieve a fast convergence.
\end{abstract}	


	\section{Introduction}
	
	Many real world problems, such as reservoir simulations, involve multiple scales and high contrast. In order to recover all the details of the media properties, one needs to
	adopt a very fine grid, which will inevitably lead to large  dimensional linear system that is hard or even impossible to solve. In order to alleviate the computational
	burden, researchers developed a lot of model reduction approaches, such as
	upscaling and multiscale methods. For example, in upscaling methods \cite{durfolsky1991homo,wu2002analysis}, one homogenizes the media properties based
	on some rules and then solve the problems on a coarse grid.
	In multiscale methods \cite{egw10,efendiev2009multiscale,Arbogast_two_scale_04,chung2015mixed,chen2003mixed,jennylt03,Wheeler_mortar_MS_12, Arbogast_PWY_07}, one still solves the problems on a coarse grid but with
	precomputed multiscale basis functions that carry small scale information of the
	media.
	
	In this paper, we  present an enrichment algorithm  in the framework of mortar mixed finite element method in solving  flow problems in heterogeneous media. We first
	compute a coarse grid solution with offline basis, which are chosen as standard polynomials basis functions. Then we iteratively compute basis
	functions based on the previous solution in the online stage, thus we call it an online method.
	The method in this paper is an extension of the online Generalized Multiscale Finite Element Methods (GMsFEM) \cite{online_cg,online_mixed,online_dg} to the mortar mixed case.
	
	Mortar mixed finite element methods \cite{mixed_dd1988,arbogast2000} are a modification of mixed finite element methods by introducing a Lagrange multiplier to impose the continuity of flux. These methods enjoy some advantages, such as  mass conservation which is very important in flow problems,  nonconforming grid discretization, and allowing domain decomposition setting which  yields a symmetric and positive definite bilinear form that defined only on the
	interfaces of the grid.  Recently, multiscale mortar mixed finite element methods \cite{Arbogast_PWY_07,arbogast2013ms} were designed to
	reduce the degree of freedom of the mortar mixed finite element method and provide
	an approximate solution. In these methods, the construction of the mortar space is a key part. Polynomial or homogenized multiscale  basis functions
	are used to form the mortar space. However, polynomials are only sufficient for very smooth
	media,  while homogenized multiscale functions lack global information which is still insufficient for accurate simulation
	of coupled flow and transport  problems(the transport velocity is a solution to a heterogeneous flow problem) if there are long range channels inside
	the media \cite{xiao2013multiscale}. Developing  efficient  domain decomposition preconditioners \cite{arbogast2013ms} to solve the fine scale problem is a good method to tackle
	this kind of problem. However, if one can equip the mortar space
	with basis functions  that can capture global information with the evolving of time, then we can avoid solving  the fine problems.
	Basis functions  with limited global information \cite{efendiev2006accurate} had proven to be
	an effective strategy using mixed finite element and finite volume coupling for the flow and transport problem. However, it fails to work in the framework of mortar mixed coupling based on our numerical study.
	Online basis functions  developed in this paper are  able to keep global information of complicated media,
	and we can use multiple online basis functions  to compute the velocity for the transport equation.
	
	Residual driven GMsFEM is an iterative algorithm that drives coarse-grid solutions
converging to the fine-grid solutions, see \cite{online_cg,online_mixed,online_dg} for the case of  finite element coupling method. The algorithm essentially includes the following steps: (1) compute an initial coarse grid solution with offline
	basis which can be polynomials or multiscale basis functions, (2) for each edge-wised coarse
	neighborhood, compute online basis by solving a homogeneous Dirichlet problem with local residual as source,
(3) compute new solution with updated basis space and then return to step (2) until
	a residual is less than the use-defined threshold.
	We  also propose an oversampling online algorithm motivated by the
	restricted domain decomposition \cite{cai1999restricted}. The idea
	is quite similar to the offline oversampling \cite{hw97}, that is
	we use a larger domain than an edge-based neighborhood to solve the local problem,
	and then take the restriction of the solution on the coarse edge as
	online basis function. This small modification turns out to be very effective in
	terms of iteration number since it includes distant information, and removes
	some boundary effects.

 We present some numerical results to show the convergence behavior of
 the method for various heterogeneous permeability fields. We  study the influence of
 the local problem size to the convergence speed.
 We  also investigate  the effects of different number of initial basis functions and  different order of contrast of the media.
We apply our approach to solve two-phase flow and transport problems.
 In two-phase flow and transport, we solve the transport equation with finite volume method on a fine grid. We adopt the online algorithm to compute the basis functions  of the initial model, and use these basis functions  to solve the flow equation on a coarse grid without adding online basis in time. We show that by adding a small number of online basis functions, the coarse-grid solution can approximate the fine-scale solution very well.

The paper is organized as follows. In section \ref{pre}, we first describe the coarse and fine discretization of the domain, then present the framework of mortar mixed finite element method, followed by the description of the domain decomposition method. In section \ref{online basis}, we introduce the iterative algorithm together with some analysis. In section \ref{oversampling}, we present an oversampling online method. Numerical examples are given in section \ref{numerical-results}, and conclusions are made in the last section.

\section{Preliminaries}\label{pre}
We consider the following second order elliptic equation in mixed formulation:
\begin{subequations}\label{problem}
	\begin{alignat}{2}
	\label{original equation-1}
	\boldsymbol{u} + \kappa\nabla p &= 0 \qquad && \text{in $\Omega$,}\\
	\label{original equation-2}
	\nabla \cdot \boldsymbol{u}  &= f && \text{in $\Omega$,}\\
	\label{boundary condition}
	\boldsymbol{u\cdot n} &= 0 && \text{on $\partial \Omega$,}
	\end{alignat}
\end{subequations}
where $\Omega \subset \mathbb{R}^d (d=2, 3)$ is a bounded polyhedral domain with outward unit normal vector $\boldsymbol{n}$ on the boundary,
$f \in L^2(\Omega)$, $\kappa$ represents the permeability  field that may vary over multiple spacial scales.

\subsection{Coarse and fine grids}
The online basis functions are constructed locally on a coarse grid. In this section we introduce coarse and fine grids. Let $\Omega$ be divided into non-overlapping polygonal coarse blocks
$K_i$ with diameter $H_i$ so that $\overline{\Omega}=\cup_{i=1}^N\overline{K}_i$,
where $N$ is the number of coarse blocks.
The decomposition of the domain can be nonconforming. We call $E_H$ a coarse
edge of the coarse block $K_i$ if $E_H = \partial K_i \cap \partial K_j $ or $E_H= \partial K_i \cap \partial{\Omega}$.
Let $\mathcal{E}_H(K_i)$ be the set of all coarse edges on the boundary of the coarse block $K_i$ and $\mathcal{E}_H=\cup_{i=1}^N\mathcal{E}_H(K_i)$
be the set of all coarse edges.

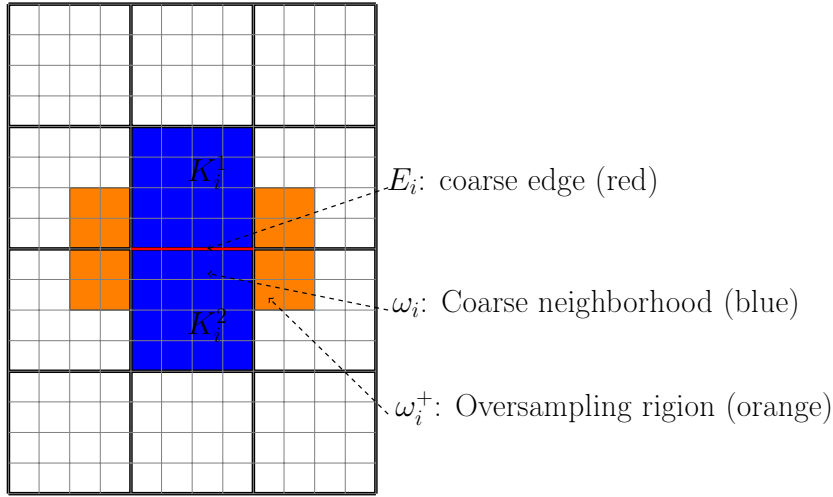
\begin{figure}[H]
	\centering
	\resizebox{0.7\textwidth}{!}{
		\begin{tikzpicture}[scale=0.7]
		\filldraw[fill=orange, draw=black] (2,6) rectangle (10,10);
		\filldraw[fill=blue, draw=black] (4,4) rectangle (8,12);
		
		\draw[step=4,black, line width=0.8mm] (0,0) grid (12,16);
		\draw[step=1,gray, thin] (0,0) grid (12,16);
		\draw[ultra thick, red](4, 8) -- (8,8);
		\draw [->,dashed, thick](12.4,10) -- (6.5, 8);
		\node at (16.8,10.2)  {\huge $E_{i}$: coarse edge (red)};
		
		\draw [->,dashed, thick](12.4,6) -- (6.5, 7.2);
		\node at (19.2,6.2)  {\huge $\omega_{i}$: Coarse neighborhood (blue)};
		
		\draw [->,dashed, thick](12.4,2.6) -- (8.5, 6.4);
		\node at (19.8,2.8)  {\huge $\omega_{i}^+$: Oversampling rigion (orange)};
		\node at (6.5,10.5)  { \huge $K_i^1$} ;
		\node at (6.5,5.5)  { \huge $K_i^2$} ;
		\end{tikzpicture}
	}
	\caption{ Illustration of a coarse edge $E_i$, and its coarse neighborhood $\omega_{i}$, oversampling rigion $\omega_{i}^+$.}
	\label{grids}
\end{figure}

We further partition each each coarse block $K_i$ into a finer mesh with mesh size $h_i$.
Let $\mathcal{T}_h=\cup_{i=1}^N\mathcal{T}_h(K_i)$ be the union of all these partitions, which is a fine mesh partition of the domain $\Omega$.
We use $h=\text{max}_{1\leq i\leq n}h_i$ to denote the mesh size of $\mathcal{T}_h$.
In addition, we let
$\mathcal{E}_h(K_i)$
be the set of all edges of the partition $\mathcal{T}_h(K_i)$ and $\mathcal{E}_h^0(K_i)$ be the set of all interior edges of the partition
$\mathcal{T}_h(K_i)$ and let $\mathcal{E}_h=\cup_{i=1}^N\mathcal{E}_h(K_i)$ be the set of all edges in the partition $\mathcal{T}_h$. 
Figure~\ref{grids} gives an illustration of the constructions of the two grids. The black lines represent the coarse grid, and the gray lines represent the fine grid.
For each coarse edge $E_i$, we define a coarse neighborhood $\omega_i$ as the union of all coarse blocks having the edge $E_i$.
Figure~\ref{grids} shows a coarse neighborhood $\omega_i$ in the blue  color.

%

\subsection{Variational form}
We introduce the following spaces
$$L_2(K_i)= \left\{p: \int_{K_i} p^2 <\infty \right\},$$
$$H(\text{div}; K_i)= \left\{ {\bf v}\in L_2(K_i)^d: \text{div} ({\bf v})\in L_2(K_i)\right\}.$$
Denote $(\cdot, \cdot)_{K_i}$ for the $L_2(K_i)$ or $L_2(K_i)^d$ inner product, and $\left\langle \cdot, \cdot\right\rangle_{\partial K_i}$ for the duality pairing on boundaries and interfaces, $d$ is the dimension of the space.
For each subdomain $i$, define
$${\bf V}_i= \left\{ {\bf v} \in H(\text{div}; K_i): {\bf v}\cdot {\bf n}|_{\partial \Omega\cap \partial K_i}=0 \right \} \quad \text{and} \quad {\bf V}=\oplus_{i=1}^N {\bf V}_i,$$
$$W_i=L_2(K_i) \quad \text{and} \quad W= \left\{ w\in L_2(\Omega): \int_{\Omega}w =0\right\},$$
$$M_{i}=H^{1/2}(E_{i}), \quad \text{and} \quad M= \oplus_{i=1}^N M_{i}.$$

The variational form for the system (\ref{original equation-1})-(\ref{boundary condition})  using mortar mixed finite element method is formulated as: find
${\bf u}\in {\bf V}$, $p\in W$ and $\lambda \in M$ such that for each $1\leq i \leq N$,
\begin{subequations}\label{Eq: variation form}
\begin{alignat}{2}
\label{variation-1}
(\kappa^{-1} \boldsymbol{u}, \boldsymbol{v})_{K_i} - (p, \nabla\cdot\boldsymbol{v})_{K_i}+\left\langle\lambda, \boldsymbol{v}\cdot \boldsymbol{n}_i\right\rangle_{ E_i} &=0\quad && \forall ~ \boldsymbol{v}\in \boldsymbol{V}_i,\\
\label{variation-2}
(\nabla\cdot \boldsymbol{u}, w)_{K_i}&= (f,w)_{K_i}\quad&& \forall ~ w\in  W_i,\\
\label{variation-3}
\sum_{i=1}^N\left\langle \boldsymbol{u}\cdot \boldsymbol{n}_i, \mu\right\rangle_{E_i} &= 0\quad && \forall ~ \mu\in M.
\end{alignat}
\end{subequations}

\subsection{The finite element approximation}
Let ${\bf V}_{h,i}\times W_{h,i}\subset {\bf V}_{i}\times W_{i}$ be any of the mixed finite element spaces satisfying the inf-sup condition for which $\nabla\cdot {\bf V}_{h,i}= W_{h,i}$, e.g., the Raviart-Thomas spaces. Define
 ${\bf V}_h=\oplus_{i=1}^N {\bf V}_{h,i}$ and $W_h=\oplus_{i=1}^N
 W_{h,i}/\mathbb{R}$ for the global discrete flux and pressure.  Let $M_{H,i}, M_{h,i}\subset L_2(E_{i})$ be the local coarse and fine mortar finite space respectively, and
$M_H=\oplus_{1\leq i \leq N}M_{H,i}, M_h=\oplus_{1\leq i \leq N}M_{h,i}$ be the entire coarse and fine mortar finite element spaces. We also denote the restriction of $M_h$ on $\mathcal{E}_H$ as $M_H^f$,
which implies $M_H$ is a subspace of $M_H^f$.

We formulate the finite element approximation as: find
${\bf u}_h\in {\bf V_h}$, $p_h\in W_h$ and $\lambda_H \in M_H$ such that for each $1\leq i \leq N$,
\begin{subequations}\label{Eq: discrete form}
\begin{alignat}{2}
\label{discrete-1}
(\kappa^{-1} \boldsymbol{u}_h, \boldsymbol{v}_h)_{K_i}-(p_h, \nabla\cdot\boldsymbol{v}_h)_{K_i} +\left\langle\lambda_H, \boldsymbol{v}_h\cdot \boldsymbol{n}_i\right\rangle_{ E_i}&=0 \quad && \forall ~\boldsymbol{v}_h\in \boldsymbol{V}_{h,i},\\
\label{discrete-2}
(\nabla\cdot \boldsymbol{u}_h, w_h)_{K_i}&= (f,w_h)_{K_i}\quad&& \forall ~ w_h\in  W_{h,i},\\
\label{discrete-3}
\sum_{i=1}^N\left\langle \boldsymbol{u}_h\cdot \boldsymbol{n}_i, \mu_H\right\rangle_{E_i} &= 0\quad && \forall ~\mu_H\in M_H.
\end{alignat}
\end{subequations}
We note that the coarse mortar space is used in this system. Similar system holds using fine mortar space $\lambda_h\in M_h.$
 Local conservation is enforced by (\ref{discrete-2}), and (\ref{discrete-3}) enforces weak continuity of flux across the interfaces with respect to the mortar space $M_H.$

\subsection{Interface problem}

The main feature of the mortar mixed finite element method is that it could be implemented by just solving a global system on the coarse mesh together with the solutions of some
local problems.

Define bilinear forms $a_{H,i}: M_{H,i}\times M_{H,i} \rightarrow \mathbb {R}, i=1,\cdots, N$ by
$$a_{H,i}(\lambda, \mu)= - \left\langle {\bf u}_h^{\ast}(\lambda)\cdot {\bf n}_i, \mu \right\rangle|_{E_i},$$ and $a_H: M_H\times M_H \rightarrow \mathbb {R}$ by
$$a_H=\sum_{i=1}^N a_{H,i}(\lambda, \mu), $$
where $\big({\bf u}_h^{\ast}(\lambda), p_h^{\ast}(\lambda)\big)\in  {\bf V}_{h} \times W_{h}$ solves ($\lambda$ given, $f=0$)
\begin{subequations}\label{Eq: domain dec1}
\begin{alignat}{2}
\label{domain dec1-1}
\big(\kappa^{-1} \boldsymbol{u}_h^{\ast}(\lambda), \boldsymbol{v}_h\big)_{K_i}-\big(p_h^{\ast}(\lambda), \nabla\cdot\boldsymbol{v}_h\big)_{K_i}&= -\left\langle\lambda, \boldsymbol{v}_h\cdot \boldsymbol{n}_i\right\rangle_{ E_i} \qquad && \forall ~ \boldsymbol{v}_h\in \boldsymbol{V}_{h,i},\\
\label{domain dec1-2}
\big(\nabla\cdot \boldsymbol{u}_h^{\ast}(\lambda), w_h\big)_{K_i}&= 0\qquad&& \forall ~ w_h\in  W_{h,i},
\end{alignat}
\end{subequations}
for each $1\leq i\leq N.$

Define linear functionals $g_{H,i}: M_{H,i} \to \mathbb {R}$ by
$$g_{H,i}(\mu)= \left\langle  \boldsymbol{\bar{u}}_h \cdot {\bf n}_i, \mu \right\rangle|_{E_i},$$ and $g_H: M_H \to \mathbb {R}$ by
$$g_H(\mu)=\sum_{i=1}^N g_{H,i}(\mu),$$
where $(\boldsymbol{\bar{u}}_h, \bar{p}_h)\in  {\bf V}_{h} \times W_{h}$ solves ($\lambda=0, f$ given) for $1\leq i\leq N$
\begin{subequations}\label{Eq: domain dec2}
\begin{alignat}{2}
\label{domain dec2-1}
\big(\kappa^{-1} \boldsymbol{\bar{u}}_h, \boldsymbol{v}_h\big)_{K_i}-\big(\bar{p}_h, \nabla\cdot\boldsymbol{v}_h\big)_{K_i}&= 0 \quad && \forall ~ \boldsymbol{v}_h\in \boldsymbol{V}_{h,i},\\
\label{domain dec2-2}
\big(\nabla\cdot \boldsymbol{\bar{u}}_h, w_h\big)_{K_i}&= (f,w_h)_{K_i}\quad&& \forall ~ w_h\in  W_{h,i}.
\end{alignat}
\end{subequations}

Define the coarse variational interface problem about the mortar pressure as: find $\lambda_H\in M_H$ such that
\begin{equation}\label{interface}
a_H(\lambda_H, \mu)=g_H(\mu) \quad \forall ~\mu \in M_H.
\end{equation}
It is proven in \cite{arbogast2000} that the interface problem (\ref{interface})  produces the solution of
(\ref{discrete-1})-(\ref{discrete-3}) via
$$\boldsymbol{u}_h= {\bf u}_h^{\ast}(\lambda)+\boldsymbol{\bar{u}}_h, p_h=\tilde{p}_h-\frac 1 {|\Omega|}\int_{\Omega}\tilde{p}_h, $$
where $\tilde{p}_h=p_h^{\ast}(\lambda)+\bar{p}_h. $

The solution of the interface problem (\ref{interface}), interpreted from the  point view of  multiscale method, is to construct multiscale basis functions over the coarse blocks.  First we design a basis for $M_H.$  For each interface $E_{i}$, from the set of mortar basis $\lambda_H$ associated with this interface, we can obtain the multiscale basis ${\bf u}_h^{\ast}(\lambda_H)$ over the coarse domains $K_{i_1}$ and $K_{i_2}.$ From these, we get a system of equations from  (\ref{interface}) directly, and solve it in any appropriate way.

The interface bilinear form $a_H(\cdot,\cdot)$ is symmetric and positive semi-definite on $M_H$ and this system can be solved by preconditioned conjugate gradient method. See \cite{Arbogast_PWY_07,hdgunified} and reference therein for more details.
The goal of this paper is to design residual driven based online enriched space $M_H$. We will use the notation $a_h(\cdot,\cdot)$ and $g_h(\cdot)$ if $a_H(\cdot,\cdot)$ and $g_H(\cdot)$ act on the space $M_h$.
\begin{remark}
	One can obtain the snapshot solution by taking $M_H=M_H^f$ and solve the problem on the coarse mesh . By solving the system corresponding to
 $a_h(\xi_h, \mu): M_h\times M_h\to R$
 equals the linear form $ g_h(\mu): M_h \to R$, we can get the fine scale solution. Note that
 the snapshot solution is equivalent to the fine scale solution.
\end{remark}

\subsection{Offline space}
To obtain online basis functions by using residuals, we first compute an offline solution from an offline mortar space $M_H^{\text{off}}$.
There are various choices for  the offline space, the simplest one is polynomials on the
coarse edge $E$. Another choice is trigonometric functions. One can also
use homogenized multiscale basis (see \cite {arbogast2015two,xiao2013multiscale}) and GMsFEM based multiscale basis (see \cite{mortaroffline}). We remark that using offline  multiscale basis is more expensive than using
continuous basis like polynomials. In this paper we consider using polynomial functions as offline basis functions to compute the initial offline
solution for simplicity and for cheap computational cost. In the next section, we discuss the construction of online basis by using residuals.

\section{Residual driven online basis}\label{online basis}
	
Using the offline space mentioned earlier is a promising choice in various scenarios.
However, in some applications such as the reservoir simulation, it is very hard to obtain a
satisfiable velocity field for the accurate simulation of the transport of flows in highly heterogeneous media with only a small number of offline basis  (see \cite {arbogast2015two,xiao2013multiscale}). Even the multiscale basis may fail to work in some cases since the offline basis only includes local information of the media.
Therefore it is important to
construct multiple multiscale basis functions  that can capture the global information of the complicated media systemically. Next, we describe our algorithm to construct such locally-supported basis functions  in the coarse grid. These new basis functions are computed in the online stage of computations. Therefore, we call them online basis functions. Next, we describe the algorithm for constructing these online basis.

\subsection{Algorithm}
Before we present the algorithm, we first introduce some notations. Denote the seminorm $||\lambda||_{M_{H}^f}=a_H^f(\lambda,\lambda)^{1/2}$ induced by $a_H^f(\cdot,\cdot)$ on $M_H^f$, we also define
linear functional on $M_H^f$ by $R_H^f(\lambda)=g_H^f(\lambda)-a_H^f(\lambda_{\text{ms}},\lambda)$.
The norm we use is given by
\begin{equation*}
||R_{\omega_i}||^2_{M_{H}^{f,*}}=\displaystyle \sup_{\lambda\in M_H^f}\frac{R_H^f(\lambda)}{||\lambda||_{M_{H}^f}}
\end{equation*}
For each coarse edge $E_i$, let $\omega_i=K_i^1\cup K_i^2$ be its neighborhood (see Figure.\ref{grids}).
Let $M_{H,i}^f$ be the restriction of $M_H^f$ on the $\cup_{i=1,2}\partial K_i$,
then we define $a_{H,i}^f(\cdot,\cdot)$ be the restriction of $a_H^f(\cdot,\cdot)$
on $M_{H,i}^f$, and $R_{H,i}^f(\lambda)$ be the restriction of $R_H^f(\lambda)$ on $M_{H,i}^f$. Similarly, the norm of this subspace is $||\lambda||_{M_{H,i}^f}=a_{H,i}^f(\lambda,\lambda)^{1/2}$. The norm for $R_{H,i}^f(\lambda)$ can be defined by
\begin{equation*}
||R_{\omega_i}||^2_{M_{H,i}^{f,*}}=\displaystyle \sup_{\lambda\in M_{H,i}^f}\frac{R_{H,i}^f(\lambda)}{||\lambda||_{M_{H,i}^f}}
\end{equation*}
We will iteratively enrich the offline space by constructing new online basis functions based on the solution computed in the previous solution space.
Let the index $l\ge0$ be the enrichment level. At the level $l$, we use $M_H^l$ to denote the corresponding
global coarse space and $\lambda_{\text{ms}}^l$ is the corresponding solution on the coarse edges.
$M_{H,i}^l$ is again the restriction of $M_H^l$ on $\cup_{i=1,2}\partial K_i.$ $M_H^0$ consists of the offline basis, i.e., polynomials, while 
the space $M_H^l (l\geq 1)$ contains both the offline and online basis functions. For each $l=0,1,2,\cdots,$ we perform
the following calculations:
\newline
\textbf{Online iterative algorithm:}
\newline
Step 1: Find the multiscale solution in the current space $M_H^l$. That is, find $\lambda_{\text{ms}}^l\in M_H^l$ such
that $a_H(\lambda_{\text{ms}}^l,\lambda)=g_H(\lambda)$ for all $\lambda\in M_H^l$.
\newline
Step 2: Pick non-overlapping neighborhoods. We select non-overlapping neighborhoods $\omega_1$, $\omega_2 $,
$\cdots$, $\omega_I$ $\subseteq \Omega$.
\newline
Step 3: Compute online basis. For each $\omega_i$, we solve for $\mu_{H,i}\in M_{H,i}^f$ such that
\begin{equation*}
a_{H,i}^f(\mu_{H,i},\lambda)=R_{H,i}^f(\lambda)\quad\quad \forall\lambda\in M_{H,i}^f,
\end{equation*}
Those $\mu_{H,i}$'s  are the new online basis. The new coarse mortar space can be updated
by setting $M_H^{l+1}=M_H^l\bigoplus \text{span}\{{\mu_{H,1},\mu_{H,2},...,
	\mu_{H,I}}\}$.

We will repeat these steps until error indicator is small or we have reached certain number of
basis functions.

	In the above algorithm, all the computations are performed on the space $M_H^f$, we do not recover the
	full solution until the last step. To achieve this, we need to assemble the corresponding
	finite element matrix of the bilinear form $a(\cdot,\cdot)$ and linear functional $g(\cdot)$ on the space $M_H^f$, denoted by $A_H^f$ and $F_H^f$ respectively.
	These involve  solving $n_{ke}^i$ zero source Dirichlet boundary value problems  for
	each coarse block $K_i$, where $n_{ke}^i$ is the number of fine grid edges on the boundary of $K_i$.
	 This may be expensive (although it can be parallelized naturally), one can also consider another
	 equivalent algorithm which will be introduced in Section \ref{oversampling}.
\begin{remark}
		 Step 3 is equivalent to solve a local zero Dirichlet boundary condition problem with the
		 full recovered local residual as source. But if we have assembled matrix $A_H^f$  and  factorized local component of $A_H^f$ for each coarse edge before the
		 online iteration. Then iterative computation cost is cheap since in this case the dominant
		 computation can be done  with parallelization before the online iteration.
\end{remark}

\begin{remark}
Here we only consider the uniform enrichment, one can also do adaptive enrichment by using an error indicator and setting a pre-defined
 tolerance to decide which coarse neighborhoods need to be enriched, see (\cite{online_mixed}).
\end{remark}

\section{Residual driven online basis with oversampling}\label{oversampling}	

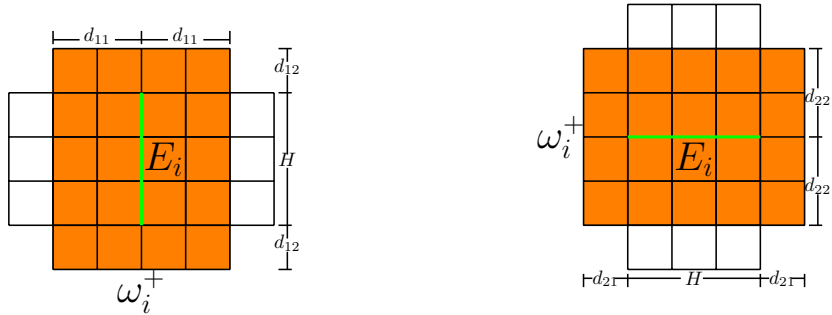
\begin{figure}[H]
		\centering
		\resizebox{0.7\textwidth}{!}{
			\begin{tikzpicture}[scale=0.8]
			\filldraw[fill=orange, draw=black] (1,0) rectangle (5,5);
			\draw[step=1,black, thick] (1,0) grid (5,5);
			\draw[step=1,black, thick] (0,1) grid (6,4);
			\draw[ultra thick, green](3, 1) -- (3,4);
			\node at (3.5,2.5)  { \huge $E_{i}$} ;
			\node at (12.5,3)  { \huge $\omega_{i}^+$} ;
			\draw[ thick, black](1, 5.1) -- (1,5.4);
			\draw[ thick, black](1, 5.2) -- (1.6,5.2);
			\node at (2,5.3)  { $d_{11}$};
			\draw[ thick, black](2.4, 5.2) -- (3,5.2);
			\draw[ thick, black](3, 5.1) -- (3,5.4);
			\draw[ thick, black](3, 5.2) -- (3.6,5.2);
			\node at (4,5.3)  { $d_{11}$};
			\draw[ thick, black](4.4, 5.2) -- (5,5.2);
			\draw[ thick, black](5, 5.1) -- (5,5.4);
			
			\draw[ thick, black](6.1, 5) -- (6.4,5);
			\draw[ thick, black](6.3, 5) -- (6.3,4.7);
			\node at (6.3,4.6)  { $d_{12}$};
			\draw[ thick, black](6.1, 4) -- (6.4,4);
			\draw[ thick, black](6.3, 4.4) -- (6.3,4);
			\draw[ thick, black](6.3, 4) -- (6.3,2.6);
			\node at (6.3,2.5)  { $H$};
			\draw[ thick, black](6.3, 2.3) -- (6.3,1);
			\draw[ thick, black](6.1, 1) -- (6.4,1);
			\draw[ thick, black](6.3, 1) -- (6.3,0.7);
			\node at (6.3,0.6)  { $d_{12}$};
			\draw[ thick, black](6.3, 0.4) -- (6.3,0);
			\draw[ thick, black](6.1, 0) -- (6.4,0);

			\filldraw[fill=orange, draw=black] (13,1) rectangle (18,5);
			\draw[step=1,black, thick] (13,1) grid (18,5);
			\draw[step=1,black, thick] (14,0) grid (17,6);
			
			\draw[ultra thick, green](14, 3) -- (17,3);
			\node at (15.5,2.5)  { \huge $E_{i}$} ;
			\node at (3,-0.5)  { \huge $\omega_{i}^+$} ;
			\draw[ thick, black](14,-0.1) -- (14,-0.4);
			\draw[ thick, black](13,-0.1) -- (13,-0.4);
			\draw[ thick, black](17,-0.1) -- (17,-0.4);
			\draw[ thick, black](18,-0.1) -- (18,-0.4);
			\draw[ thick, black](13, -0.2) -- (13.2,-0.2);
			\node at (13.5,-0.2)  { $d_{21}$};
			\draw[ thick, black](13.7, -0.2) -- (14,-0.2);
			\draw[ thick, black](14, -0.2) -- (15.2,-0.2);
			\draw[ thick, black](15.7, -0.2) -- (17,-0.2);
			\node at (15.5,-0.2)  { $H$};
			\draw[ thick, black](17, -0.2) -- (17.2,-0.2);
			\draw[ thick, black](17.7, -0.2) -- (18,-0.2);
			\node at (17.5,-0.2)  { $d_{21}$};
			
			\draw[ thick, black](18.1,1) -- (18.4,1);
			\draw[ thick, black](18.1,3) -- (18.4,3);
			\draw[ thick, black](18.1,5) -- (18.4,5);
			\draw[ thick, black](18.3,5) -- (18.3,4.1);
			\draw[ thick, black](18.3,3.7) -- (18.3,3);
			\node at (18.3,3.9)  { $d_{22}$};
			\draw[ thick, black](18.3,3) -- (18.3,2.1);
			\draw[ thick, black](18.3,1.7) -- (18.3,1);
			\node at (18.3,1.9)  { $d_{22}$};
			\end{tikzpicture}
		}
		\caption{ Illustration of an oversampled neighborhood associated with the coarse edge $E_i$.}
		\label{oversampgrid}
	\end{figure}

We can also apply the idea of oversampling \cite{hw97} in the above algorithm. The implementation
of oversampling we  introduce here is different from the no oversampling case.
We can not directly use the residual that only defined on $M_H^f$. Instead, we need to recover the global
residual that lives on $M_h$ by solving a zero source Dirichlet problem on each coarse element $K_i$.
Again, we denote the seminorm $||\lambda||_{M_{h}}=a_h(\lambda,\lambda)^{1/2}$ induced by $a_h(\cdot,\cdot)$ on $M_h$, we also define
linear functional on $M_h$ by $R_h(\lambda)=g_h(\lambda)-a_h(\lambda_{\text{ms}},\lambda)$.
We note here that the $a_h(\cdot,\cdot)$ and $g_h(\cdot)$
require solve a number of zero source non-homogeneous Dirichlet problems for each fine grid element in $\mathcal{T}_h$. However, there is some linear relationship between the solution and the coefficient (single constant), so the computation is cheap. We denote $A_e$ and $F_e$ be finite element matrix corresponding
to $a_h(\cdot,\cdot):M_h\times M_h\to R$ and $g_h(\cdot):M_h\to R$.
Our goal is still to find an online basis to enrich $M_H$ which is defined on
all coarse edges, however we will no longer solve the local problem in step 3 on
$\omega_i$, instead we will use a sightly  larger domain to compute the local online basis. More specifically, for each coarse edge $E_i$, we consider
a domain $\omega_i^+\supset E_i$ (see Figure~\ref{oversampgrid} for the illustration
of $\omega_i^{+}$) as the target local domain to perform local computation.
Let $M_h^{\omega_i^{+}}$ be the restriction of $M_h$ on $\omega_i^+$. Let $a_{h,i}(\cdot,\cdot)$ and $g_{h,i}(\cdot)$ be the restriction of $a_h(\cdot,\cdot)$ and $g_h(\cdot)$ on $M_h^{\omega_i^{+}}$.
 \newline
We keep the notation in Section \ref{online basis}, then we have
 \newline
 \textbf{Oversampling  online iterative algorithm:}
 \newline
 Step 1: Find the multiscale solution in the current space $M_H^l$. That is, find $\lambda_{\text{ms}}^l\in M_H^l$ such
 that $a(\lambda_{\text{ms}}^l,\lambda)=g(\lambda)$ for all $\lambda\in M_H^l$.
 \newline
 Step 2: Pick oversampled neighborhoods. For each coarse edge $E_i$, we select an oversampled 
 coarse neighborhood $\omega^+_i$ (see Figure \ref{oversampgrid}). We repeat this selection for 
 coarse edges $E_1, E_2, \cdots E_I$. Then we obtain oversampled neighborhoods $\omega_1^+$, $\omega_2^+ $,
 $\cdots$, $\omega_I^+$ $\subseteq \Omega$. The index $\{1,2,...,I\}$ can be chosen such that $\displaystyle\bigcup_{i=1,\cdots,I}\omega_i$ (not $\omega_i^+$)
 form a non-overlapping partition of $\Omega$.
 \newline
Step 3: Compute the global full scale solution. Compute the global solution $\lambda_{h}(\xi) \in M_{h}$  by solving  zero source problem
with the restriction of $\lambda_{\text{ms}}^l$ on coarse blocks as Dirichlet boundary conditions.
 \newline
 Step 4: Compute online basis. For each $\omega_i^+$, we solve for $\mu_{h,i}^+\in M_{h,i}$ such that
 \begin{equation*}
 a_{h,i}(\mu_{h,i}^+,\lambda)=R_{h,i}(\lambda)\quad\quad \forall\lambda\in M_{h,i}.
 \end{equation*}
  \newline
 Step 5: Take the restriction. We take the restriction of $\mu_{h,i}^+$ on coarse edges $E_i$, denoted
 by $\mu_{H,i}$

Those $\mu_{H,i}$'s  are the new oversampling online basis, and then the new basis space can be updated accordingly by adding them to the previous solution space.
We can also pre-compute and factorize the matrix associated with  $a_{h,i}(\cdot,\cdot)$, whose
computational cost may be cheaper than the no oversampling case, since the number of fine scale edges in $\omega_i^+$ may be less than those of $\omega$.
The major differences of oversampling and no oversampling approaches are: (1) the domain that used to compute the local
online basis for the oversampling is larger than the standard domain $\omega_i$ in terms of the direction of
$E_i$. (2) the computation of no oversampling can be done on space $M_H^f$, there is no necessary to
compute the residual on $M_h$, therefore the online iterative computation of the no oversampling case is cheaper than the 
oversampling case.

\section{Numerical examples}\label{numerical-results}
\begin{figure}[H]
	\centering
	\subfigure[$\kappa_{1}$]{
		\includegraphics[width=2.5in]{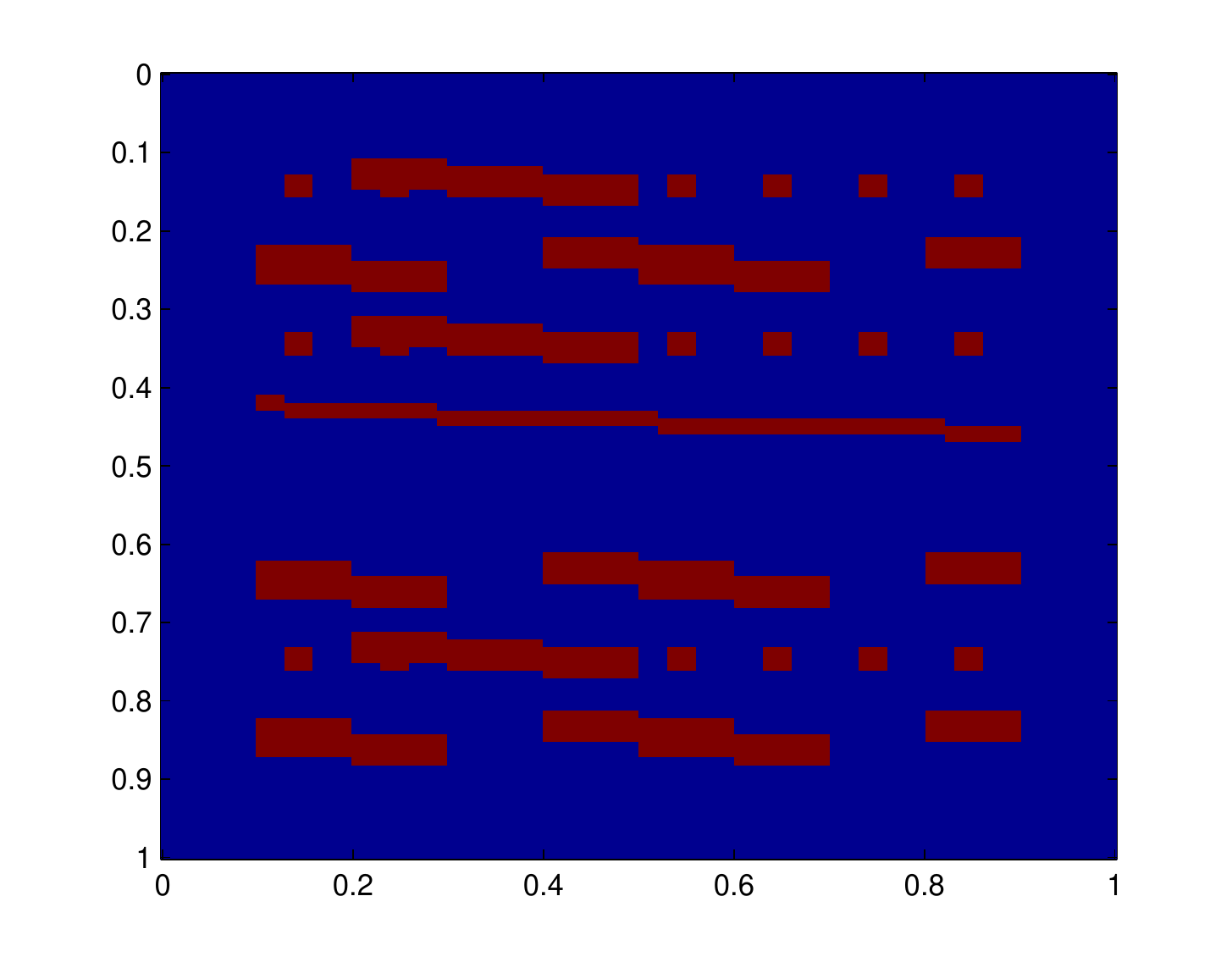}}
	\subfigure[$\kappa_{2}$ in $\log_{10}$ scale]{
		\includegraphics[width=3.5in]{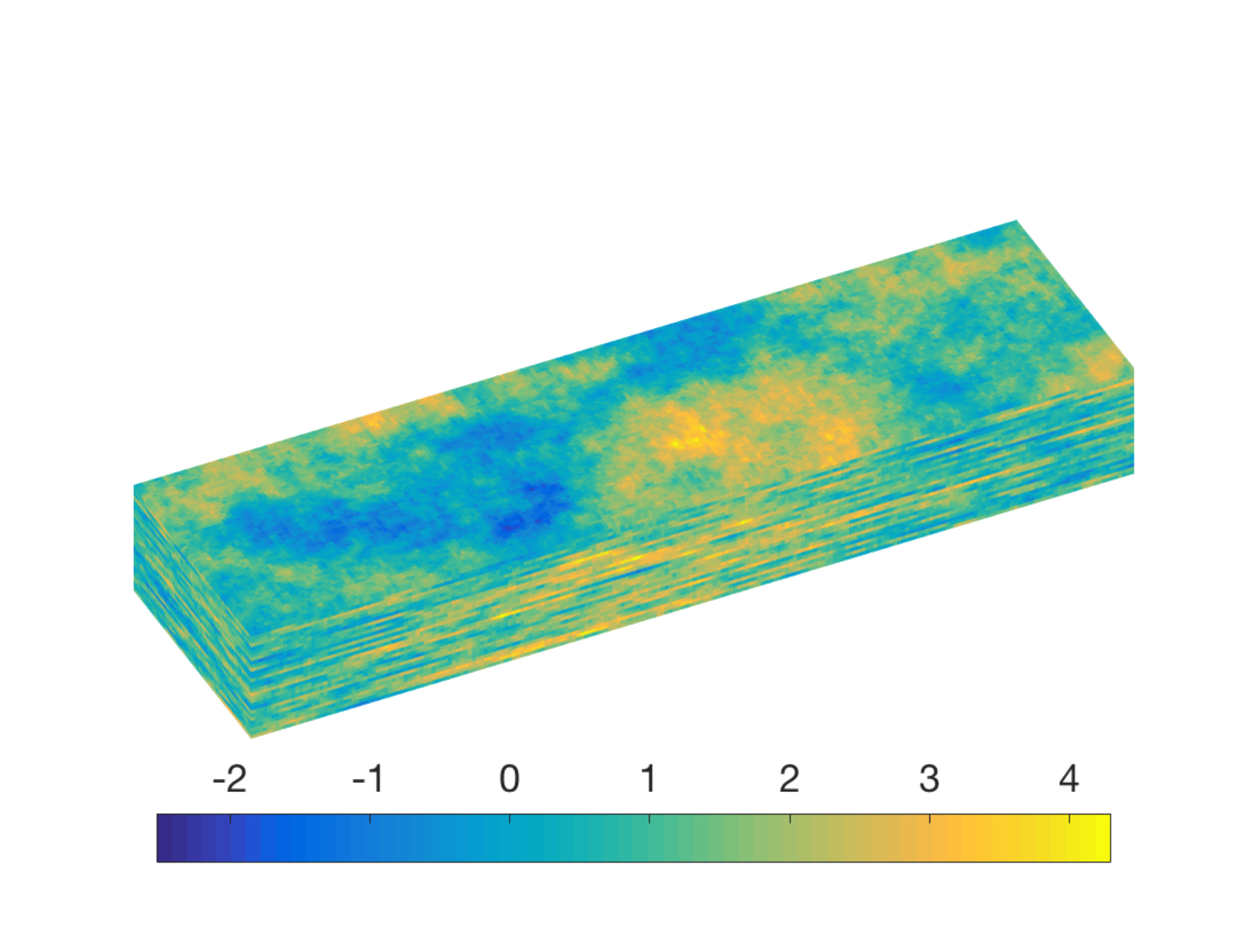}}
		\subfigure[$\kappa_{3}$ in $\log_{10}$ scale]{
			\includegraphics[width=3.5in]{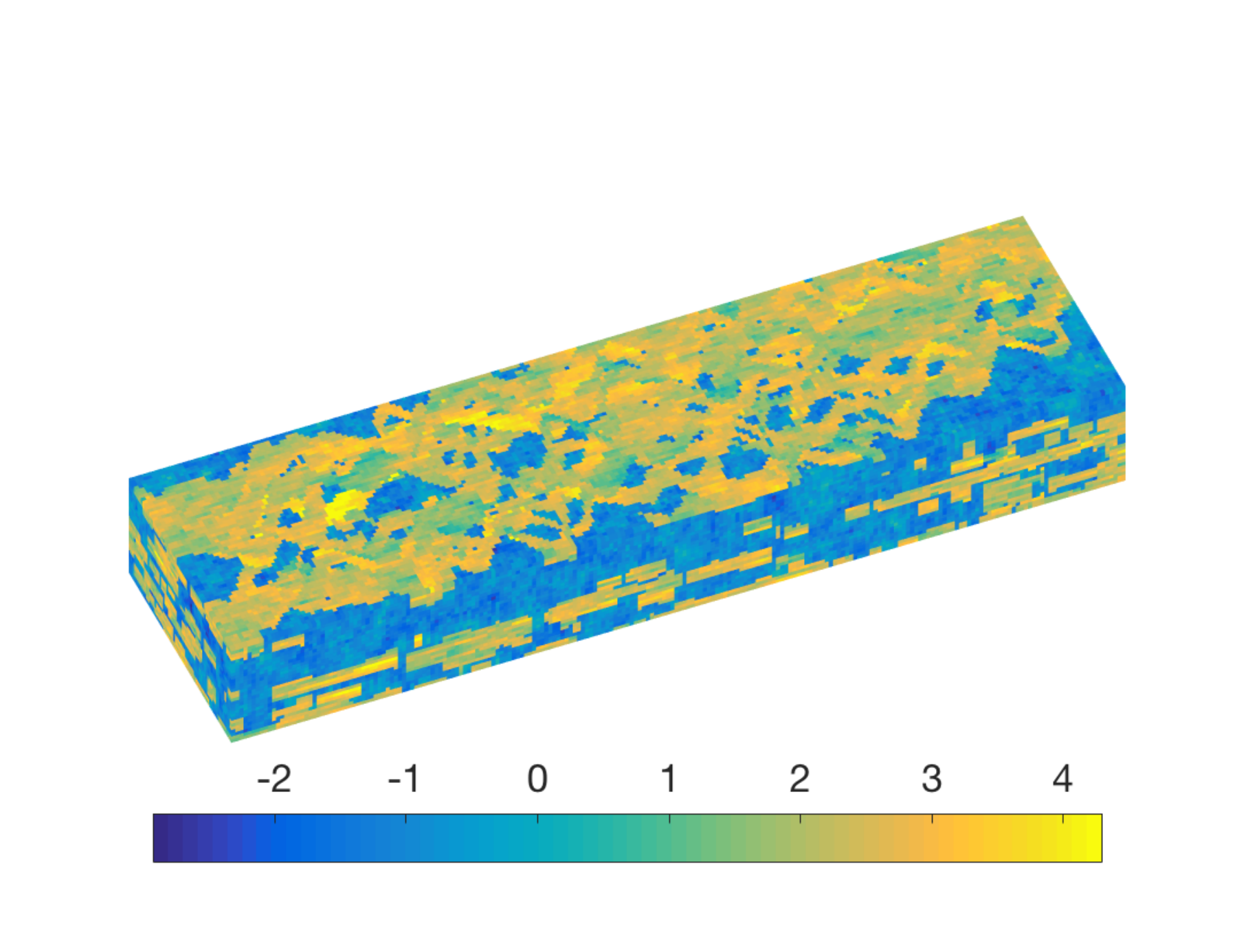}}
	\caption{Permeability fields.}
	\label{model} 
\end{figure}

In this section, we present several representative examples to show the performance of our method.
We consider three models with permeability $\kappa$ depicted in Figure \ref{model}.   For model 1 in Figure \ref{model}(a), we note that $\kappa=1$ in the blue region and $\kappa=\eta$  in the red region, where $\eta$ will be specified in each example.
As it is shown, the first model contains high contrast, long channels, and isolated inclusions. The second (the first 30  layers of the SPE10 model) and third model  (the last 30 layers of the SPE10 model) are selected from the tenth SPE comparative solution project (SPE 10)  \cite{Aarnes2005257}. 
The SPE 10 model (its full model has $60 \times 220 \times 85$ cells) is  used as a  benchmark  to test different upscaling techniques and multiscale methods,  and is therefore a good test case for our methodology.

We define the following errors for both pressure and flux to quantify the accuracy of the online
multiscale solution
\begin{equation*}
e_p:=\frac{\|p_{\text{ms}}-p_f\|_{L^2,\Omega}}{\|p_f\|_{L^2,\Omega}}, \quad  e_{\boldsymbol{u}}:=\frac{\|\boldsymbol{u}_{\text{ms}}-\boldsymbol{u}_f\|_{\kappa,\Omega}}{\|\boldsymbol{u}_f\|_{\kappa,\Omega}}
\end{equation*}
where $\|\boldsymbol{u}\|_{\kappa,\Omega}^2=\int_{\Omega}\kappa^{-1}\boldsymbol{u}^2$.

Our method is tested on elliptic problems in  Section \ref{num exam: elliptic}, and on two-phase flow and transport problems in Section \ref{two-phase}. In Section \ref{num exam: elliptic}, we show the performance of our method for elliptic problems. We see that adding a few number of online basis functions per edge is able to produce fast convergence speed. In particular, oversampling achieves even faster convergence.  Moreover, our method is robust in the sense that the convergence is independent of the order of contrast. In Section \ref{two-phase}, we present numerical results for a two-phase flow and transport problem. We only enrich the solution space of the initial problem, and use this initial solution space for the rest of the simulation along the time. Our numerical results show that the online basis functions produces accurate production file along the time.

\subsection{Online method for elliptic problems}\label{num exam: elliptic}

In this example, we compare the performance  of oversampling and non-oversampling, as well as the use of  different number of offline and online basis functions. We also test the robustness  of our method with respect to different contrast orders.  

In all simulations reported below, the computational domain $D$
is divided into $N_x\times N_y$($N_x\times N_y \times N_z$ for 3D) square coarse elements, and in each coarse element, we generate a uniform $n\times n$($n\times n\times n$) fine scale square elements. For model 1, a fixed fine-grid size with $200$ is employed. We use coarse grid size $10\times 10$. The source function $f$ is zero everywhere except that  it is taken as four on the top left fine grid, and  negative four on the bottom right fine grid cell.
For model 2 and model 3, the fine grid is of size $60\times 220\times 30$ (the dimension of the fine system is 1209600) and is divided into $6\times 22\times 3$ coarse elements. The fine-grid solution is used as the reference solution in all numerical examples.

We use $\left(\begin{array}{cc}
d_{11} & d_{12} \\
d_{21}  & d_{22} \\
 \end{array}
 \right)$ to define  the local domain (see Figure~\ref{oversampgrid}) for the computation of the online basis functions. In total, three cases  given below are considered( these three cases can be extended to 3D easily). Here as define earlier,  $n$ is the number of fine elements in a coarse block for each direction:\\
Case 1: no oversampling: $\left(\begin{array}{cc}
$n$ & 0 \\
0 & $n$\\
\end{array}
\right)\\$
Case 2: oversampling case a: $\left(\begin{array}{cc}
$n$ & 1 \\
1& $n$\\
\end{array}
\right)\\$\\
Case 3: oversampling case b: $\left(\begin{array}{cc}
$ [n/2] $ & 1 \\
1 & $ [n/2] $\\
\end{array}
\right)\\$\\

In the first case, the local domain is exactly the coarse neighborhood for a coarse edge, which is the no oversampling case. In Case 2, one layer of fine cells is added to the coarse neighborhood in the direction of the edge. Therefore, the oversampling domain is larger than the coarse neighborhood. In Case 3, the layers of cells on both sides of the coarse edge are reduced to about a half of $n$, while one layer of fine cells is added to the coarse neighborhood in the direction of this edge. Thus, the oversampling domain in this case is smaller than the coarse neighborhood.

We compare the errors  of using the three  domain cases for online basis function computation, to see the performance of oversampling and non-oversampling. All the three permeability fields $\kappa_{1}, \kappa_{2}$ and $\kappa_{3}$ are considered. The results  are given in Tables \ref{model1_1}-\ref{model3_2}.  Tables \ref{model1_1}-\ref{model1_2} are errors  for model 1 by using 1 and 2 offline polynomial basis functions, respectively.   First, we check the effect of the three computational domains. Take Table \ref{model1_1} for example. In the first column, $N_b$ stands for the number of basis functions for each subdomain, and Dof stands for degree of freedoms.  We start with 1 polynomial basis function for each local subdomain. Then we use the online method to add basis functions iteratively until the number of basis functions reaches to 7. The rest of the columns are errors from using  the  three computational domains. By comparing Case 1 and Case 2, we see that the errors decay faster for Case 2 even that only one layer of fine grid cell is added on each side in the coarse edge direction.  The domain in Case 3 is smaller than that of Case 1. However, we still get faster convergence rate from Case 3 since one layer of fine grid cell is added on each side in the coarse edge direction. By comparing Case 1 with the rest two cases, we observe that the oversampling technique generally improves the convergence.  Similar results can be seen in Table \ref{model1_2}, which are obtained by starting with 2 polynomial basis functions for each local subdomain.  Tables \ref{model2_1}-\ref{model2_2} presents errors  for model 2 by starting with 1 and 4 offline polynomial basis functions, respectively. We can get the same conclusion as for model 1.  Tables \ref{model3_1}-\ref{model3_2} presents errors  for model 3 by starting with 1 and 4 offline polynomial basis functions, respectively. Compared with model 2, the errors decay slower for model 3, since the permeability field for model 3 is a spaghetti of channelling system which is much more complicated.

Next, we check the effect of using different number of offline basis functions by looking at corresponding columns in Tables \ref{model1_1}-\ref{model1_2}. Look at the column for Case 1 in Table \ref{model1_1} and Table \ref{model1_2}, the row with $N_b=6$. In total, 6 basis functions are used in each subdomain. Therefore, the sizes of the final system are the same. However, for Table \ref{model1_1}, the 6 basis functions consists of 1 offline and 5 online, while  for 
Table \ref{model1_2}, the 6 basis functions consists of 2 offline and 4 online. The former needs one more iteration on each subdomain. For both models, we see that the online solution converges to the fine grid solution whether we start with 1 or 2 offline basis functions.

\begin{figure}[H]
	\centering
	\subfigure[ 1 offline basis function]{
		\includegraphics[width=3.0in]{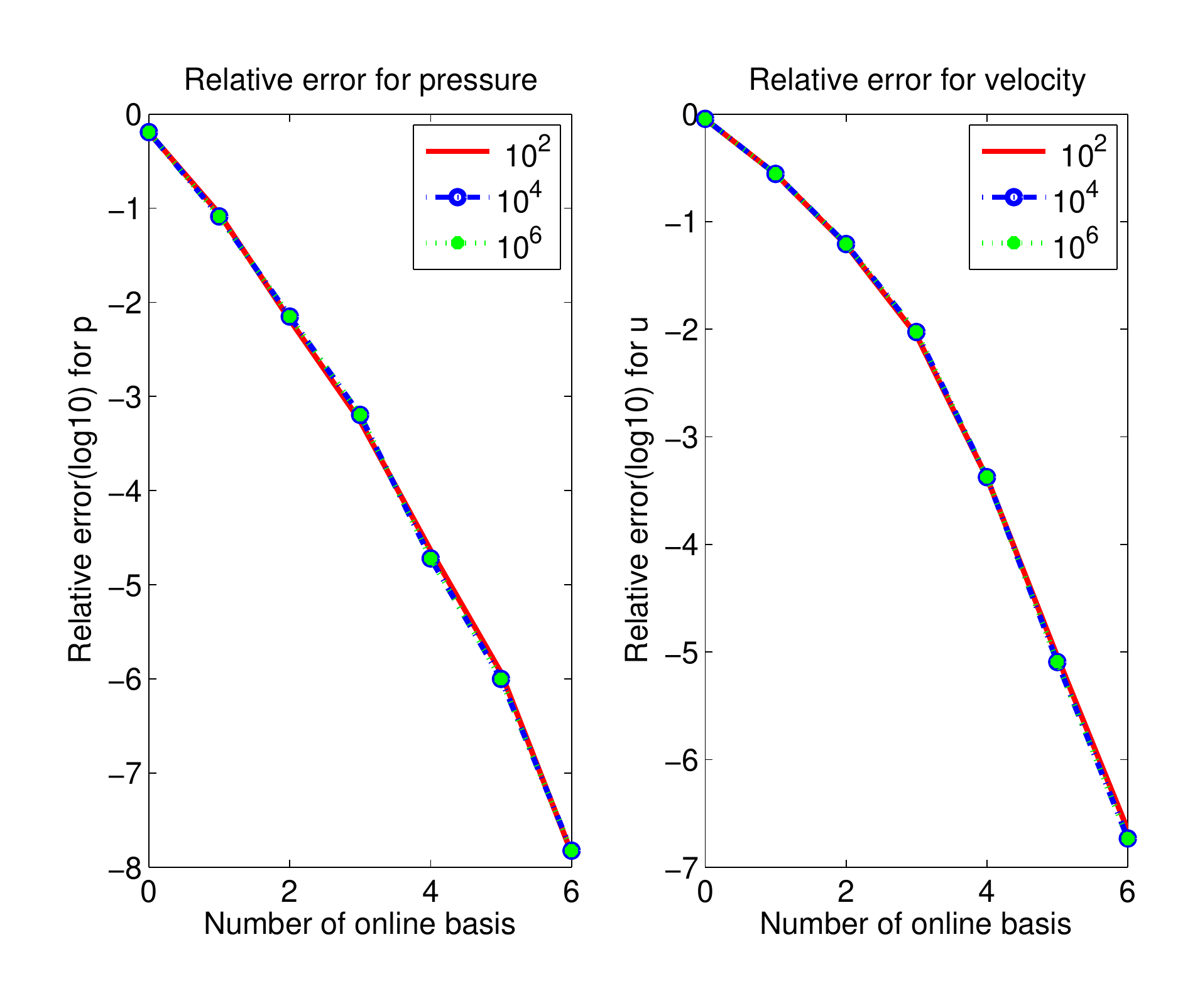}}\label{contrast_err_1_model1}
	\subfigure[2 offline  basis functions]{
		\includegraphics[width=3.0in]{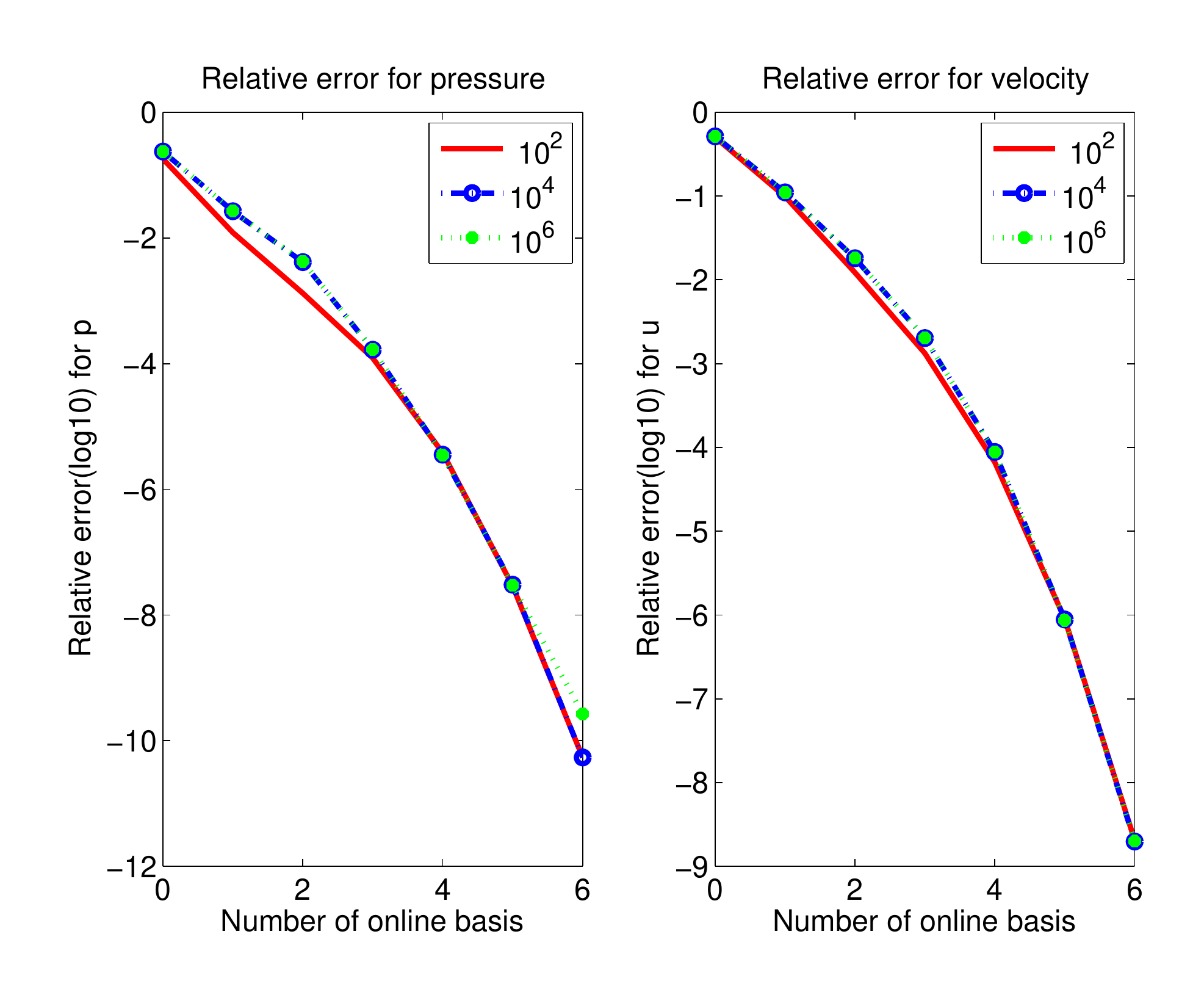}}
	\caption{Relative error of online method using different number of offline  basis functions and contrast orders $\eta=10^2, 10^4, 10^6$ for model 1.}
	\label{model1_errors} 
\end{figure}

\begin{figure}[H]
	\centering
	\subfigure[ 1 offline  basis function]{
		\includegraphics[width=3.0in]{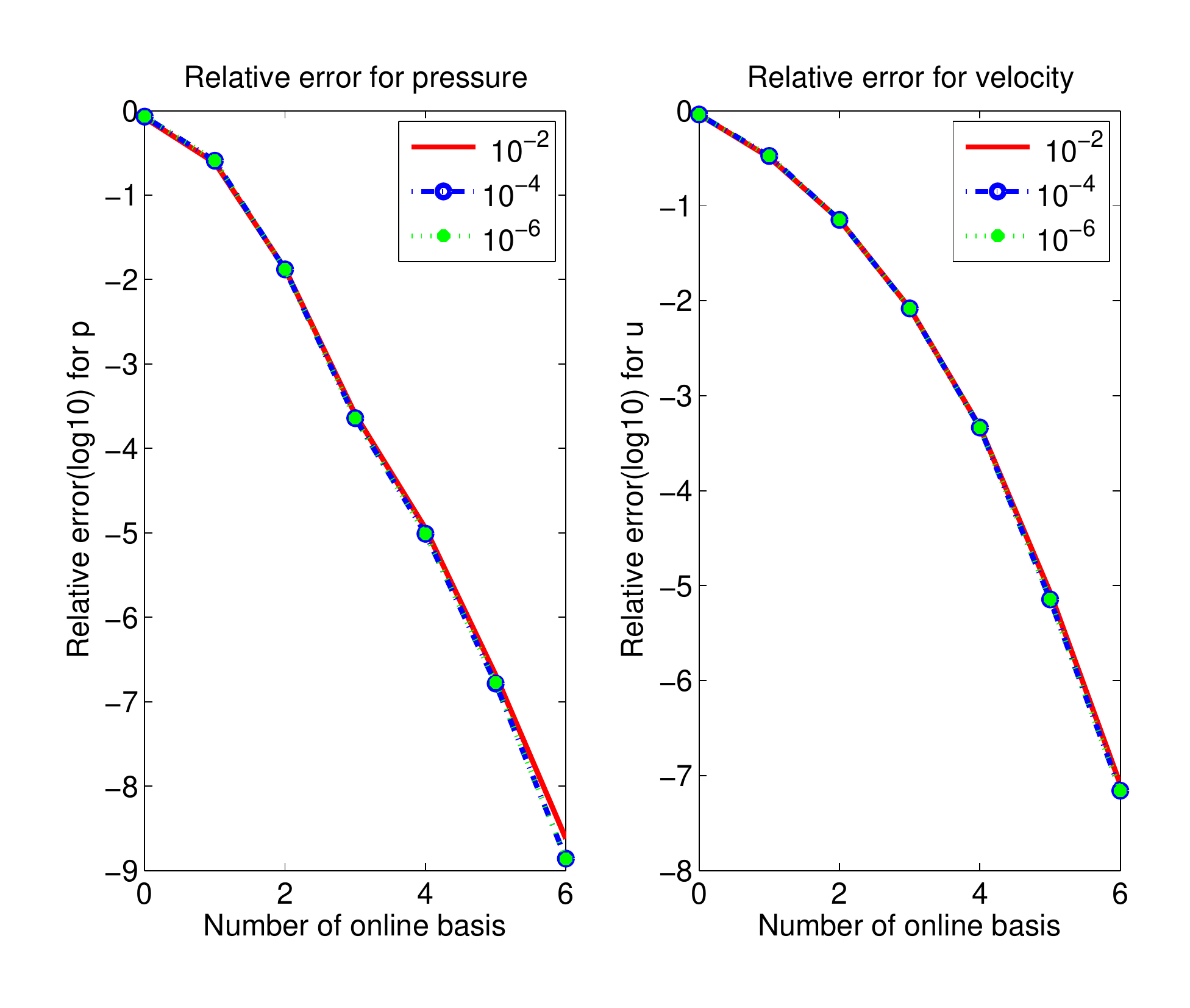}}
	\subfigure[2 offline  basis functions]{
		\includegraphics[width=3.0in]{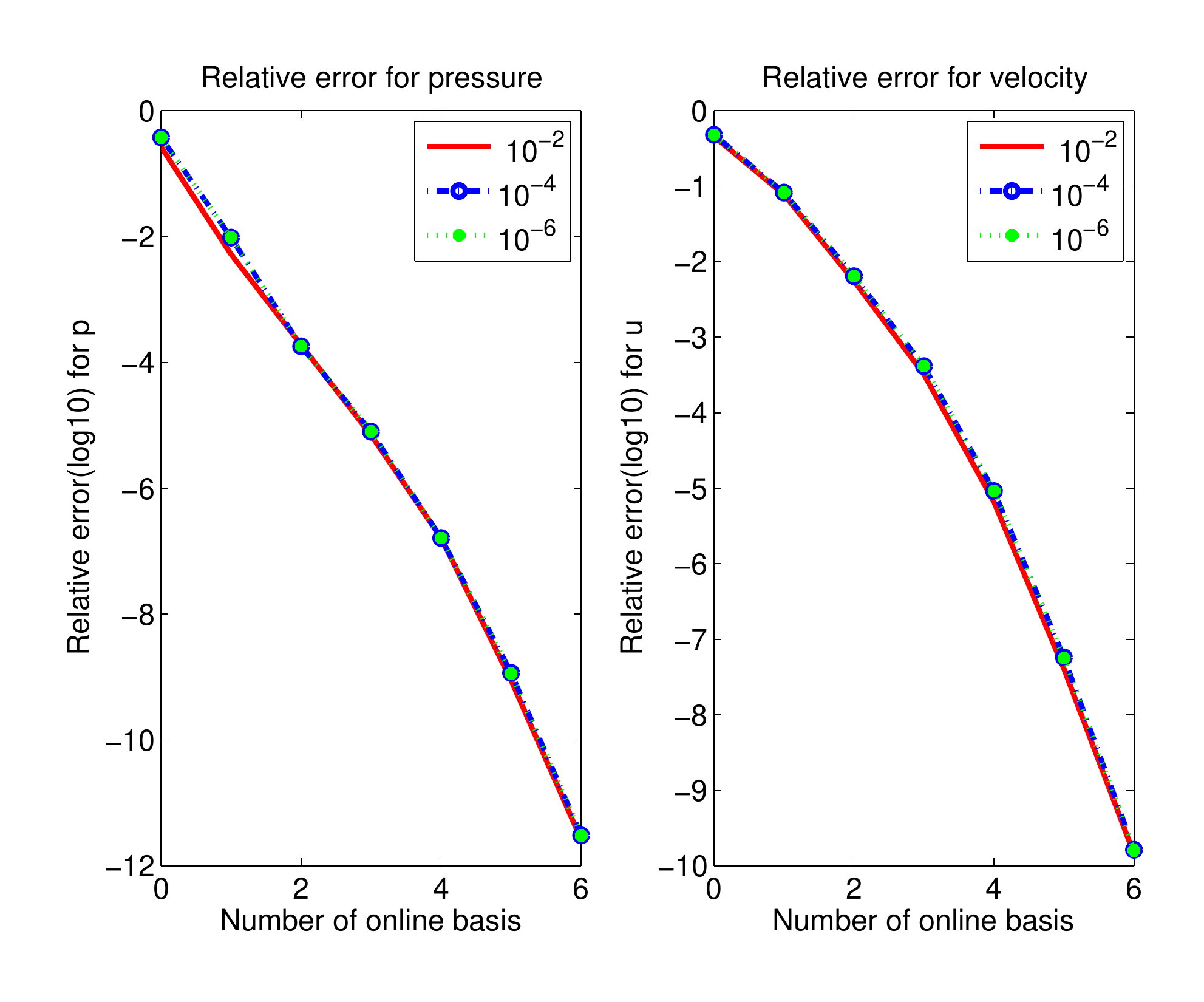}}
	\caption{Relative error of online method using different number of offline  basis functions and contrast orders $\eta=10^{-2}, 10^{-4}, 10^{-6}$ for model 1.}
	\label{model2_errors} 
\end{figure}

Finally, we show the performance for different contrast orders. We vary the order of contrast, one case is from $10^2, 10^4$ to $ 10^6$, and  the other is from $10^{-2}, 10^{-4}$ to $10^{-6}$. We plot both the pressure and velocity errors against the number of online basis functions for model 1 by using  the different  contrast values (Figures \ref{model1_errors}-\ref{model2_errors}). In Figures \ref{model1_errors} (a), we present the pressure and velocity errors against the number of online basis functions for contrast order $10^2, 10^4$ and $ 10^6$, starting with 1 offline polynomial basis. The left figure is for pressure and the right figure is for the velocity. We see that the convergence rate is almost the same for the three contrast order examples. Figure \ref{model1_errors} (b) is for the case of using 2 offline basis functions. We observe similar results, that is, the change in the contrast has almost no effect on the errors. Figure \ref{model2_errors} presents the results for contrast order $10^{-2}, 10^{-4}$ and $10^{-6}$, which also shows that the convergence lines of the three contrast cases agree well for both pressure and velocity. We conclude that the online method is robust in the sense that its convergence rate is independent of the contrast order.

\begin{table}[H]
	\centering
	 \begin{adjustbox}{max width=\textwidth}
	\begin{tabular}{|c|c|c|c|c|c|c|c|c|c|c|}
		\hline
		\multirow{2}{*}{$N_b$ (Dof)} & \multicolumn{2}{c|}{\specialcell{Case 1}} & \multicolumn{2}{c|}{ \specialcell{Case 2}} & \multicolumn{2}{c|}{\specialcell{Case 3}}\tabularnewline
		\cline{2-7}
		 & $e_p$ &$e_{\boldsymbol{u}}$ & $e_p$ &$e_{\boldsymbol{u}}$ & $e_p$ & $e_{\boldsymbol{u}}$ \tabularnewline
		\hline
		1(180)&6.42e-01&9.01e-01&6.42e-01&9.01e-01&6.42e-01&9.01e-01 \tabularnewline
		\hline
		2(360)&8.17e-02&2.78e-01 &6.32e-02&2.26e-01&6.00e-02&2.24e-01  \tabularnewline
		\hline
		3(540)&7.07e-03&6.20e-02&3.70e-03&3.44e-02&3.34e-03 &3.41e-02 \tabularnewline
		\hline
		4(720)&6.37e-04&9.42e-03&2.20e-04&2.44e-03&1.63e-04&2.00e-03 \tabularnewline
		\hline
		5(900)&1.90e-05&4.22e-04&6.06e-06&6.70e-05&3.65e-06&4.50e-05 \tabularnewline
		\hline
		6(1080)&1.00e-06&8.08e-06&2.86e-08 &3.89e-07&1.82e-08&2.64e-07  \tabularnewline
		\hline
		7(1260)&1.52e-08&1.85e-07&1.16e-11&3.41e-10&1.31e-11&4.05e-10\tabularnewline
		\hline

	\end{tabular}
\end{adjustbox}
	\caption{Relative error between multiscale solution and fine scale solution with different type of local problem cases for model 1, $N=10$,  $\eta=10^4$, using 1 offline basis function. "$N_b$" represents the total number of basis functions per coarse edge, which includes 1 offline polynomial basis function. "Dof" denotes the degree of freedom of the coarse system. }
	\label{model1_1}
\end{table}

\begin{table}[H]
	\centering
	\begin{adjustbox}{max width=\textwidth}
		\begin{tabular}{|c|c|c|c|c|c|c|c|c|c|c|}
			\hline
			\multirow{2}{*}{$N_b$ (Dof)} & \multicolumn{2}{c|}{\specialcell{Case 1}} & \multicolumn{2}{c|}{ \specialcell{Case 2}} & \multicolumn{2}{c|}{\specialcell{Case 3}}\tabularnewline
			\cline{2-7}
			& $e_p$ &$e_{\boldsymbol{u}}$ & $e_p$ &$e_{\boldsymbol{u}}$ & $e_p$ & $e_{\boldsymbol{u}}$ \tabularnewline
			\hline
		2(360) &2.37e-01 &5.15e-01 &2.37e-01 &5.15e-01&2.37e-01 &5.15e-01 \tabularnewline
		\hline
		3(540)&2.66e-02&1.11e-01&1.83e-02&9.33e-02&1.78e-02&9.00e-02 \tabularnewline
		\hline
		4(720)&4.12e-03&1.81e-02&9.77e-04&6.40e-03 &9.26e-04&6.31e-03 \tabularnewline
		\hline
		5(900)&1.67e-04&2.01e-03& 1.89e-05 &2.89e-04 &1.73e-05&2.64e-04 \tabularnewline
		\hline
		6(1080)& 3.56e-06 &8.83e-05 &6.53e-08&1.62e-06&5.84e-08&1.47e-06   \tabularnewline
		\hline
		7(1260)&3.04e-08& 8.83e-07 &1.04e-10 &2.99e-09&1.18e-10&3.14e-09 \tabularnewline
		\hline
		8(1440)&5.42e-11 &1.98e-09 &2.71e-12&1.77e-12&2.63e-12 &1.56e-12  \tabularnewline
		\hline
	\end{tabular}
\end{adjustbox}
	\caption{Relative error between multiscale solution and fine scale solution with different type of local problem cases for model 1, $N=10$, $\eta=10^4$, using 2 offline basis functions. "$N_b$" represents the  total number of basis functions per coarse edge, which includes 2 offline polynomial basis functions. "Dof" denotes the degree of freedom of the coarse system. }
	\label{model1_2}

\end{table}

\begin{table}[H]
	\centering
	\begin{adjustbox}{max width=\textwidth}
		\begin{tabular}{|c|c|c|c|c|c|c|c|c|c|c|}
			\hline
			\multirow{2}{*}{$N_b$ (Dof)} & \multicolumn{2}{c|}{\specialcell{Case 1}} & \multicolumn{2}{c|}{ \specialcell{Case 2}} & \multicolumn{2}{c|}{\specialcell{Case 3}}\tabularnewline
			\cline{2-7}
			& $e_p$ &$e_{\boldsymbol{u}}$ & $e_p$ &$e_{\boldsymbol{u}}$ & $e_p$ & $e_{\boldsymbol{u}}$ \tabularnewline
			\hline
		1(1404)&7.00e-01&7.58e-01&7.00e-01&7.58e-01&7.00e-01&7.58e-01 \tabularnewline
		\hline
		3(4212) &8.54e-02&1.55e-01&6.95e-02&1.20e-01&7.06e-02 &1.21e-01 \tabularnewline
		\hline
		5(7020)&2.66e-02&4.73e-02& 1.64e-02&2.94e-02&1.72e-02&3.01e-02 \tabularnewline
		\hline
		7(9828)&7.28e-03&1.82e-02&2.36e-03&7.94e-03&2.63e-03&8.50e-03 \tabularnewline
		\hline
		9(12636)&1.45e-03 &6.51e-03&1.28e-04& 2.00e-03 &1.46e-04&1.11e-03 \tabularnewline
		\hline
		11(15444)&7.77e-05&1.04e-03&2.51e-06 &2.68e-05 &2.28e-06&4.07e-05  \tabularnewline
		\hline
		13(18252)&8.42e-06 &1.06e-04&3.76e-08 &3.25e-07&4.49e-08 &3.49e-07  \tabularnewline
		\hline
	\end{tabular}
\end{adjustbox}
	\caption{Relative error between multiscale solution and fine scale solution with different type of local problem cases for model 2,  using 1 offline basis function. "$N_b$" represents the total number of basis functions per coarse edge, which includes 1 offline polynomial basis function. "Dof" denotes the degree of freedom of the coarse system. }
	\label{model2_1}
\end{table}

\begin{table}[H]
	\centering
	\begin{adjustbox}{max width=\textwidth}
		\begin{tabular}{|c|c|c|c|c|c|c|c|c|c|c|}
			\hline
			\multirow{2}{*}{$N_b$ (Dof) } & \multicolumn{2}{c|}{\specialcell{Case 1}} & \multicolumn{2}{c|}{ \specialcell{Case 2}} & \multicolumn{2}{c|}{\specialcell{Case 3}}\tabularnewline
			\cline{2-7}
			& $e_p$ &$e_{\boldsymbol{u}}$ & $e_p$ &$e_{\boldsymbol{u}}$ & $e_p$ & $e_{\boldsymbol{u}}$ \tabularnewline
			\hline
		4(5616) &3.04e-01&5.89e-01&3.04e-01 &5.89e-01&3.04e-01 &5.89e-01 \tabularnewline
		\hline
		6( 8424)&1.56e-02&5.98e-02&9.47e-03&4.06e-02 &9.68e-03 &4.11e-02 \tabularnewline
		\hline
		8(11232)&2.10e-03&9.67e-03&1.15e-03&5.19e-03&1.24e-03&5.50e-03\tabularnewline
		\hline
		10(14040)&5.52e-04&3.07e-03&1.73e-05 &2.31e-04&1.52e-05&2.44e-04 \tabularnewline
		\hline
		12(16848)&9.53e-06 & 2.52e-04&2.98e-07&4.13e-06&3.12e-07 &5.01e-06  \tabularnewline
		\hline
		14(19656)&3.88e-07&1.29e-05&1.22e-08&5.16e-08 &1.22e-08&5.82e-08  \tabularnewline
		\hline

	\end{tabular}
\end{adjustbox}
	\caption{Relative error between multiscale solution and fine scale solution with different type of local problem cases for model 2,  using 4 offline basis functions. "$N_b$" represents the total number of basis functions per coarse edge, which includes 4 offline polynomial basis functions.   "Dof" denotes the degree of freedom of the coarse system. }
	\label{model2_2}
\end{table}

\begin{table}[H]
	\centering
	\begin{adjustbox}{max width=\textwidth}
		\begin{tabular}{|c|c|c|c|c|c|c|c|c|c|c|}
			\hline
			\multirow{2}{*}{$N_b$ (Dof)} & \multicolumn{2}{c|}{\specialcell{Case 1}} & \multicolumn{2}{c|}{ \specialcell{Case 2}} & \multicolumn{2}{c|}{\specialcell{Case 3}}\tabularnewline
			\cline{2-7}
			& $e_p$ &$e_{\boldsymbol{u}}$ & $e_p$ &$e_{\boldsymbol{u}}$ & $e_p$ & $e_{\boldsymbol{u}}$ \tabularnewline
			\hline
		1(1404)&8.69e-01&1.15e+00&8.69e-01&1.15e+00&8.69e-01&1.15e+00\tabularnewline
		\hline
		3(4212) &5.07e-01&4.80e-01&5.09e-01&6.59e-01&5.10e-01&6.55e-01 \tabularnewline
		\hline
		5(7020)&4.62e-01&8.72e-01& 3.73e-01&3.61e-01&3.74e-01&3.60e-01 \tabularnewline
		\hline
		7(9828)&2.59e-01&3.03e-01&1.02e-01&1.72e-01&1.01e-01&1.73e-01\tabularnewline
		\hline
		9(12636)&2.33e-02 &8.98e-02&7.29e-03& 3.67e-02 &6.96e-03&3.59e-02\tabularnewline
		\hline
		11(15444)&8.13e-03&2.50e-02&3.54e-04 &3.54e-03&3.52e-04&3.49e-03 \tabularnewline
		\hline
		13(18252)&7.56e-04&7.18e-03&1.87e-05 &1.63e-04&1.62e-05 &1.80e-04  \tabularnewline
		\hline
	\end{tabular}
\end{adjustbox}
	\caption{Relative error between multiscale solution and fine scale solution with different type of local problem cases for model 3,  using 1 offline basis function. "$N_b$" represents the total number of basis functions per coarse edge, which includes 1 offline polynomial basis function. "Dof" denotes the degree of freedom of the coarse system. }
	\label{model3_1}
\end{table}

\begin{table}[H]
	\centering
	\begin{adjustbox}{max width=\textwidth}
		\begin{tabular}{|c|c|c|c|c|c|c|c|c|c|c|}
			\hline
			\multirow{2}{*}{$N_b$ (Dof)} & \multicolumn{2}{c|}{\specialcell{Case 1}} & \multicolumn{2}{c|}{ \specialcell{Case 2}} & \multicolumn{2}{c|}{\specialcell{Case 3}}\tabularnewline
			\cline{2-7}
			& $e_p$ &$e_{\boldsymbol{u}}$ & $e_p$ &$e_{\boldsymbol{u}}$ & $e_p$ & $e_{\boldsymbol{u}}$ \tabularnewline
			\hline
		4(5616) &8.36e-01&1.29e+00&8.36e-01 &1.29e+00&8.36e-01 &1.29e+00 \tabularnewline
		\hline
		6( 8424)&3.84e-01&6.28e-01&4.19e-01&7.06e-01 &4.19e-01&7.06e-01 \tabularnewline
		\hline
		8(11232)&1.87e-01&2.80e-01&9.95e-02&1.60e-01&9.96e-02&1.60e-01\tabularnewline
		\hline
		10(14040)&2.69e-02&6.89e-02&9.21e-04 &8.77e-03&9.22e-04&8.89e-03\tabularnewline
		\hline
		12(16848)&4.39e-03 & 1.17e-02&3.66e-05&4.93e-04&3.68e-05 &4.84e-04 \tabularnewline
		\hline
		14(19656)&8.82e-05&6.23e-04&4.10e-07&5.36e-06&4.54e-07&5.90e-06 \tabularnewline
		\hline

	\end{tabular}
\end{adjustbox}
	\caption{Relative error between multiscale solution and fine scale solution with different type of local problem cases for model 3,  using 4 offline basis functions. "$N_b$" represents the total number of basis functions per coarse edge, which includes 4 offline polynomial basis functions.   "Dof" denotes the degree of freedom of the coarse system. }
	\label{model3_2}
\end{table}

\subsection{A two phase flow and transport problem}\label{two-phase}
In this section, we use our method to solve a two phase flow and transport model problem. First, we summarize the underlying partial differential equations \cite{efendiev2016online,SpeJ} to simulate porous media flows. In particular, we consider two-phase flow in a reservoir domain (denoted by $\Omega$)
with the assumption that the fluid displacement is driven by viscous effects, that is, we neglect compressibility
and gravity for simplicity in our simulations. 
We consider water and oil phases which are assumed to be immiscible. By the Darcy's law, we get the following equation
for each phase

\begin{equation} \label{darcy}
{ \bf u}_l=-\frac{k_{rl}(s_l)}{\mu_l} { K} \nabla {p}
\end{equation}
where ${  \bf u}_l$ is the phase velocity, $ {K}$ is the permeability
tensor, $k_{rl}$ is the relative permeability to phase $l$ ($l=o, w$), $s_l$ is
 saturation, and ${p}$ is pressure.
Throughout the paper, we use a
single set of relative permeability.

By the mass conservation law , the following equations for the two phases are obtained:
\begin{equation} \label{mass}
\phi\frac{\partial s_l}{\partial t} + \nabla \cdot {\bf u}_l = q_l.
\end{equation}

Combining Darcy's law, mass conservation, and the property $s_w+s_o=1$, we derive the following coupled system of pressure and saturation equations (we use $s$ instead of $s_w$ for simplicity):
\begin{eqnarray}
\nabla\cdot {\bf u} &=& q_w + q_o \quad \textrm{in} \quad \Omega \label{pressEq}\\
\phi \frac{\partial s}{\partial t}+ \nabla \cdot({ f_w(s) {\bf u}}) &=& \frac{q_w}{\rho_w} \quad \textrm{in} \quad \Omega \label{satEq}\\
{\bf u}\cdot n & = & 0 \quad \textrm{on} \quad \partial {\Omega} \quad \textrm{(no flow at boundary)}\\
s(t=0) & = & s_{0} \quad \textrm{in} \quad \Omega \quad \textrm{(initial known saturation)}
\end{eqnarray}
where $\phi$ is the porosity, $\lambda$ is the total mobility defined as
\begin{equation} \label{mobility}
\lambda(s) =\lambda_w(s)+\lambda_o(s)=\frac{k_{rw}(s)}{\mu_w}+\frac{k_{ro}(s)}{\mu_o}
\end{equation}
 $f_w(s)$ is the flux function,
\begin{equation} \label{fraction}
f_w(s) =\frac{\lambda_{w}(s)}{\lambda(s)}=\frac{k_{rw}(s)}{k_{rw}(s)+\frac{\mu_w}{\mu_o}k_{ro}(s)}
\end{equation}
and ${\bf u} = { \bf u}_w + { \bf u}_o=-\lambda(s) { K}  \nabla {p}$ is the total flux. Moreover, $q_w$ and $q_o$ are volumetric source terms
for water and oil.

Here, we follow the sequential formulation, that is at each time step one solves for the pressure and velocity first and then uses the velocity to solve for the saturation. The pressure equation is solved by using the offline basis functions together with the online basis functions  computed at the initial time step,  and the saturation equation is solved by the finite volume method. We apply 
10 times cheap Jacobi \cite{mansfield1991damped} iterations to smooth the mortar multiscale solution.

The initial water saturation is taken to be zero. The velocity in Equation (\ref{satEq}) is the fine grid velocity which is obtained by projecting the multiscale velocity field  onto the fine grid. Five wells are included in the reservoir, with 1 producer in the middle and 4 injectors on the corners of the domain, i.e, the total source term is zero everywhere except that  it is taken as four on the corners of the fine grid, and  negative four in the middle fine grid cell. We report the relative saturation error at every  50 time steps, and the end of simulation time is 2500.  We define the relative saturation at time step $i$ as\begin{equation*}
e_s(i):=\frac{\|s_{\text{ms}}(i)-s_{ref}(i)\|_{L^2,\Omega}}{\|s_{ref}(i)\|_{L^2,\Omega}}.
\end{equation*}

\begin{figure}[H]
	\centering
	\subfigure[Saturation error]{
		\includegraphics[width=3.05in]{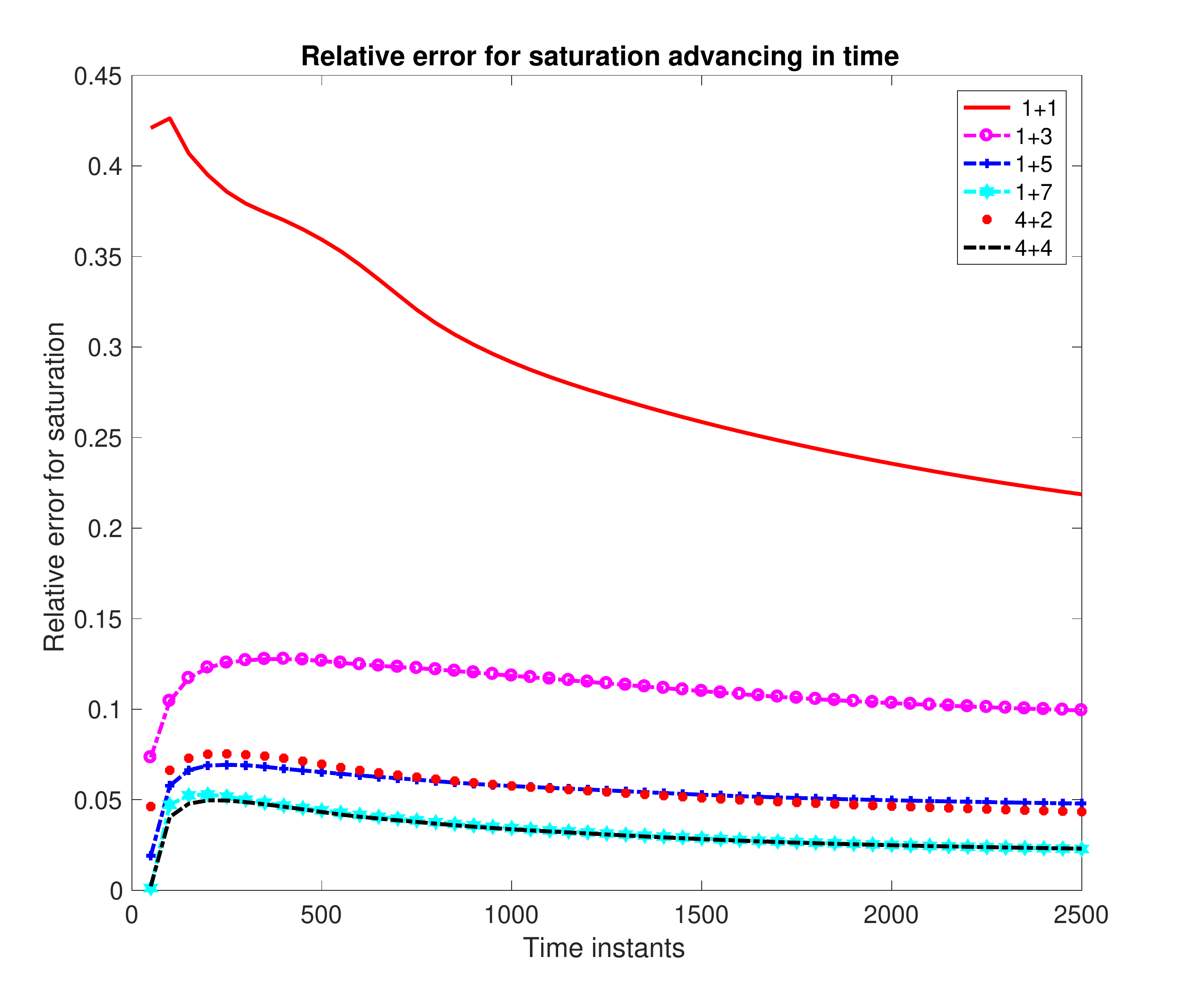}}
	\subfigure[Watercut]{
		\includegraphics[width=3.05in]{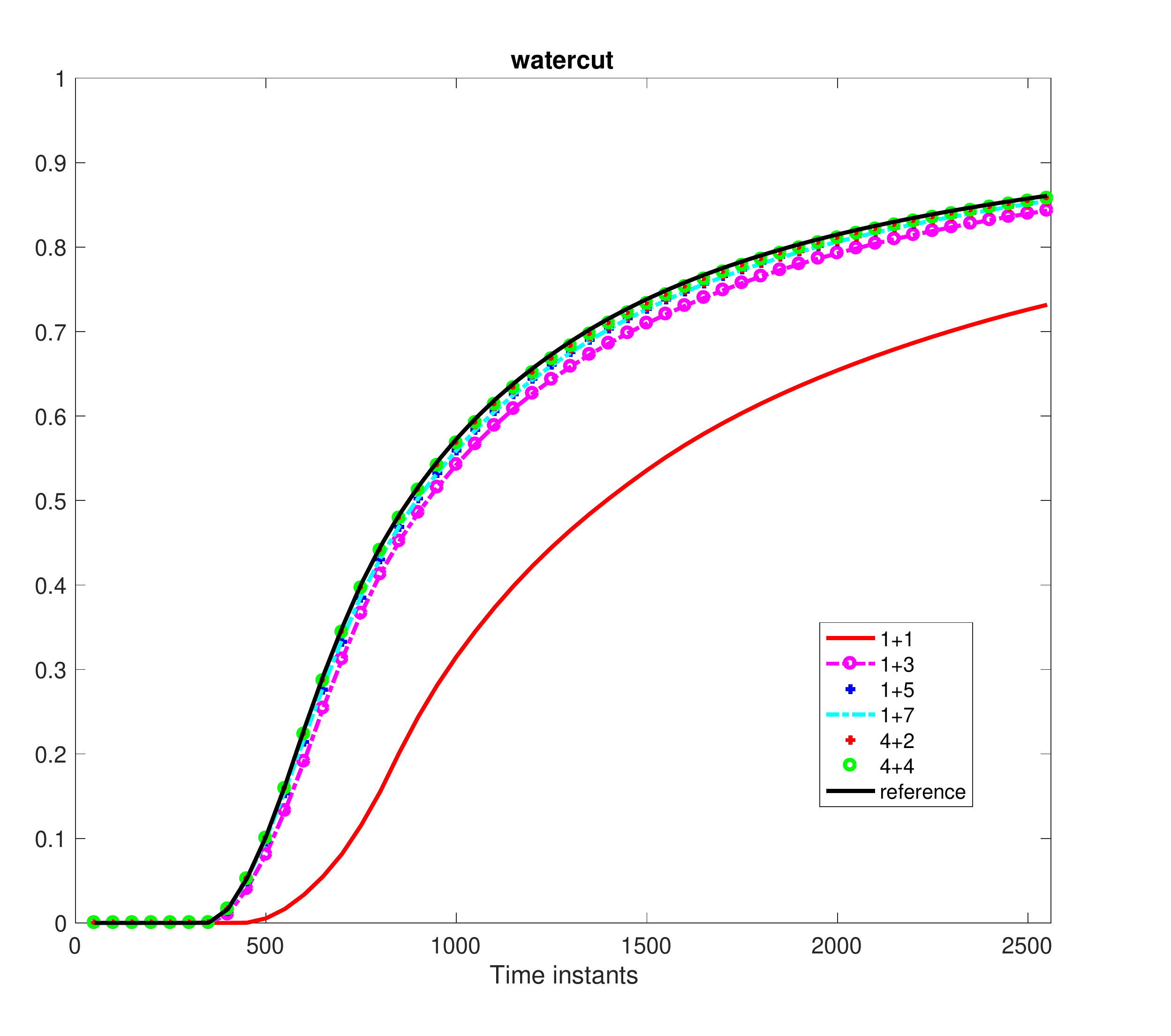}}
	\caption{(a) Saturation error advancing in time for different number offline and online basis functions; (b) Watercut for different number offline and online basis functions for the first 30 layers of SPE 10 as in Figure \ref{model}(b) }
	\label{spe_first30} 
\end{figure}

\begin{figure}[H]
	\centering
	\subfigure[ Saturation error]{
		\includegraphics[width=3.05in]{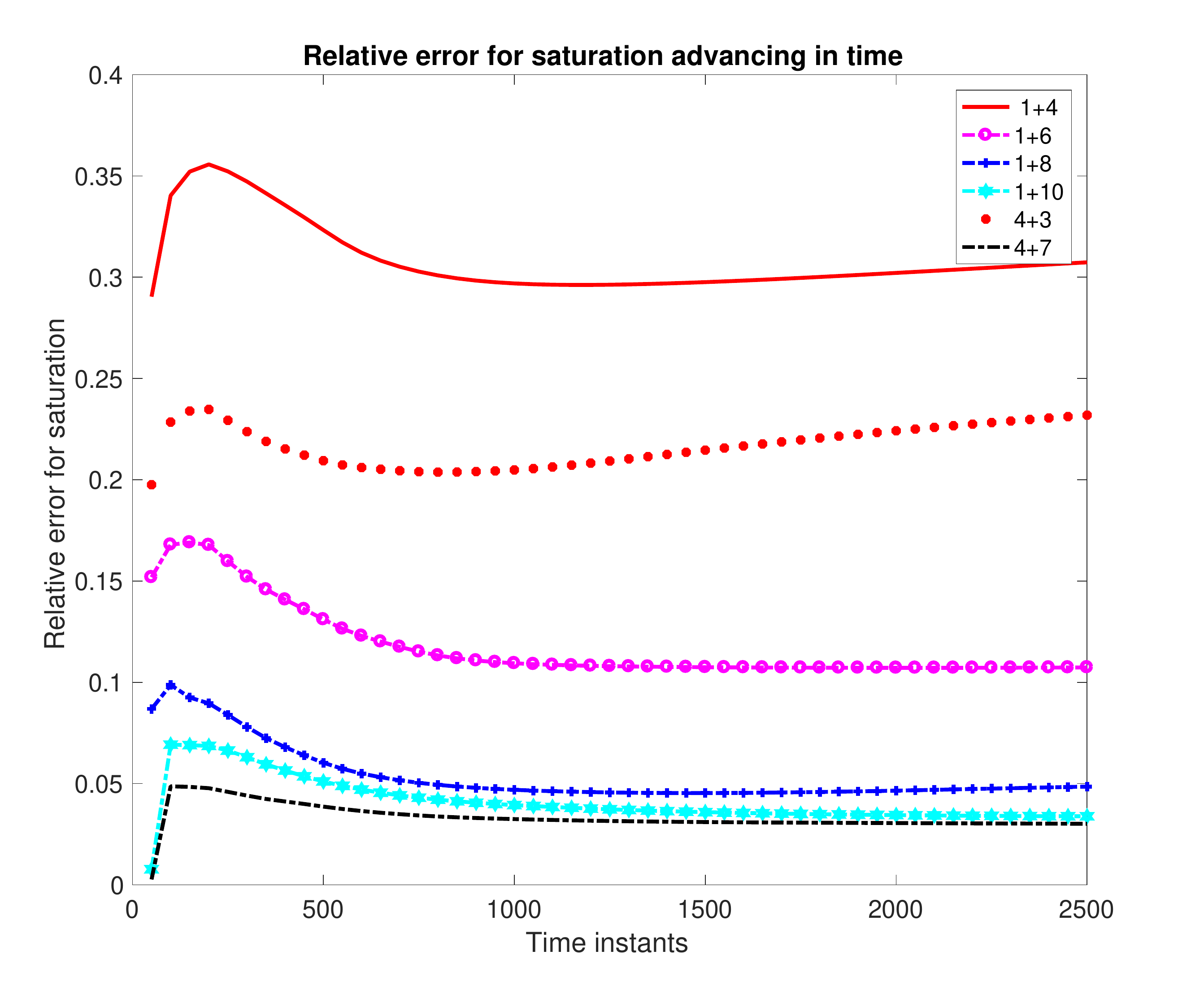}}
	\subfigure[Watercut]{
	\includegraphics[width=3.05in]{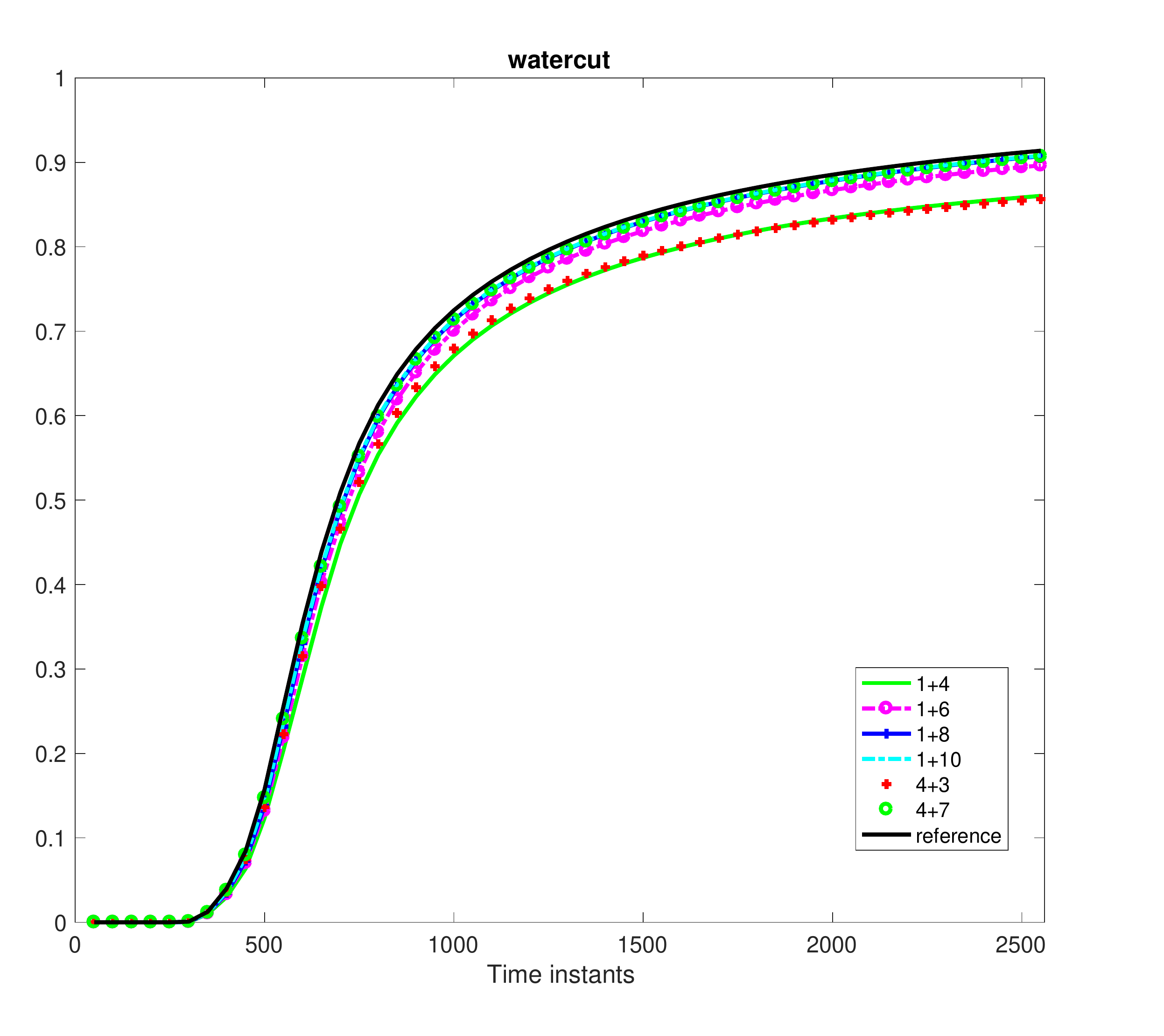}}
	\caption{(a) Saturation error advancing in time for different number offline and online basis functions; (b) Watercut for different number offline and online basis functions for the last 30 layers of SPE 10 as in Figure \ref{model}(c)}
	\label{spe_last30} 
\end{figure}

Figure \ref{spe_first30} (a) plots the relative saturation errors of adding different online basis functions to different number of offline basis function against the time instants. In the figure, in the legend the number in the form of $x+y$, the number $x$ means the number of offline basis functions, and $y$ means the number of online basis functions. For 
example, $4+2$ means 4 offline and 2 online basis functions are used. The errors are greater than $20\%$  for all the time instants if only 1 online basis is added (shown by the red line). The errors drop to less than $15\%$ if we add 3 online basis functions which is shown in the purple-circle line. The errors of the rest cases are under $10\%$.
We note that the online basis functions are only added at the initial time step, then they are fixed for the rest of the simulation time.   We also present the water-cut (water flux fractional function
$f_w(s)$) in Figure \ref{spe_first30} (b)  corresponding to the cases in Figure \ref{spe_first30} (a).  In Figure \ref{spe_first30} (b),  the red line for the case of 1 offline and 1 online basis  is far away from the black reference line, which is no surprise since we already know that the relative saturation errors are large for this case.  By adding more online basis functions, the water-cut lines get closer to the reference line.
Figure \ref{spe_last30} shows similar results for model 3, we can see that for this model, more online basis functions are needed to get satisfactory results.

\begin{figure}[H]
	\centering
	\subfigure[Reference solution ]{
		\includegraphics[width=3.0in]{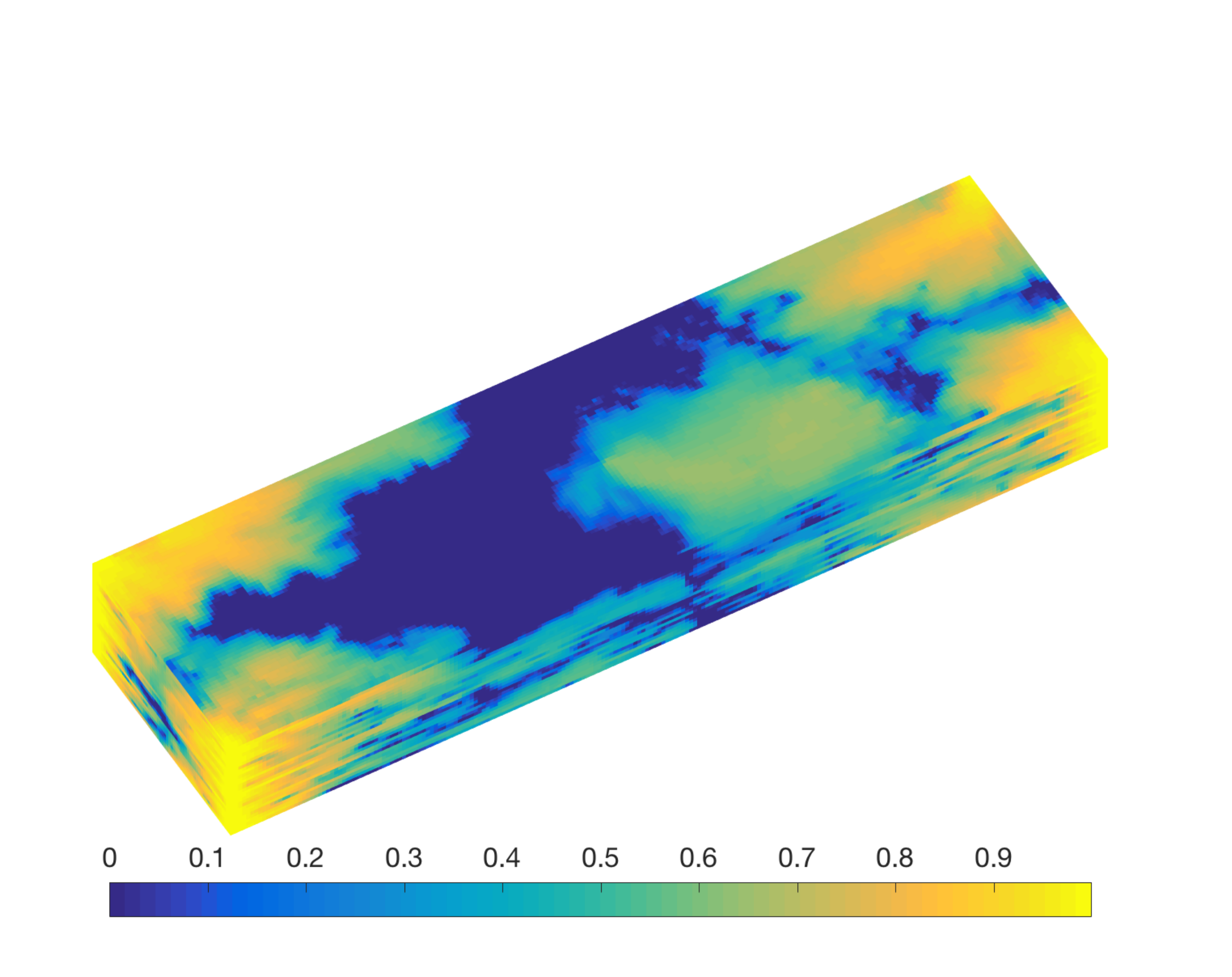}}
	\subfigure[Ms-solution with 1 offline basis and 1 online basis, relative $L^2$ error is  21.9$\%$]{
		\includegraphics[width=3.0in]{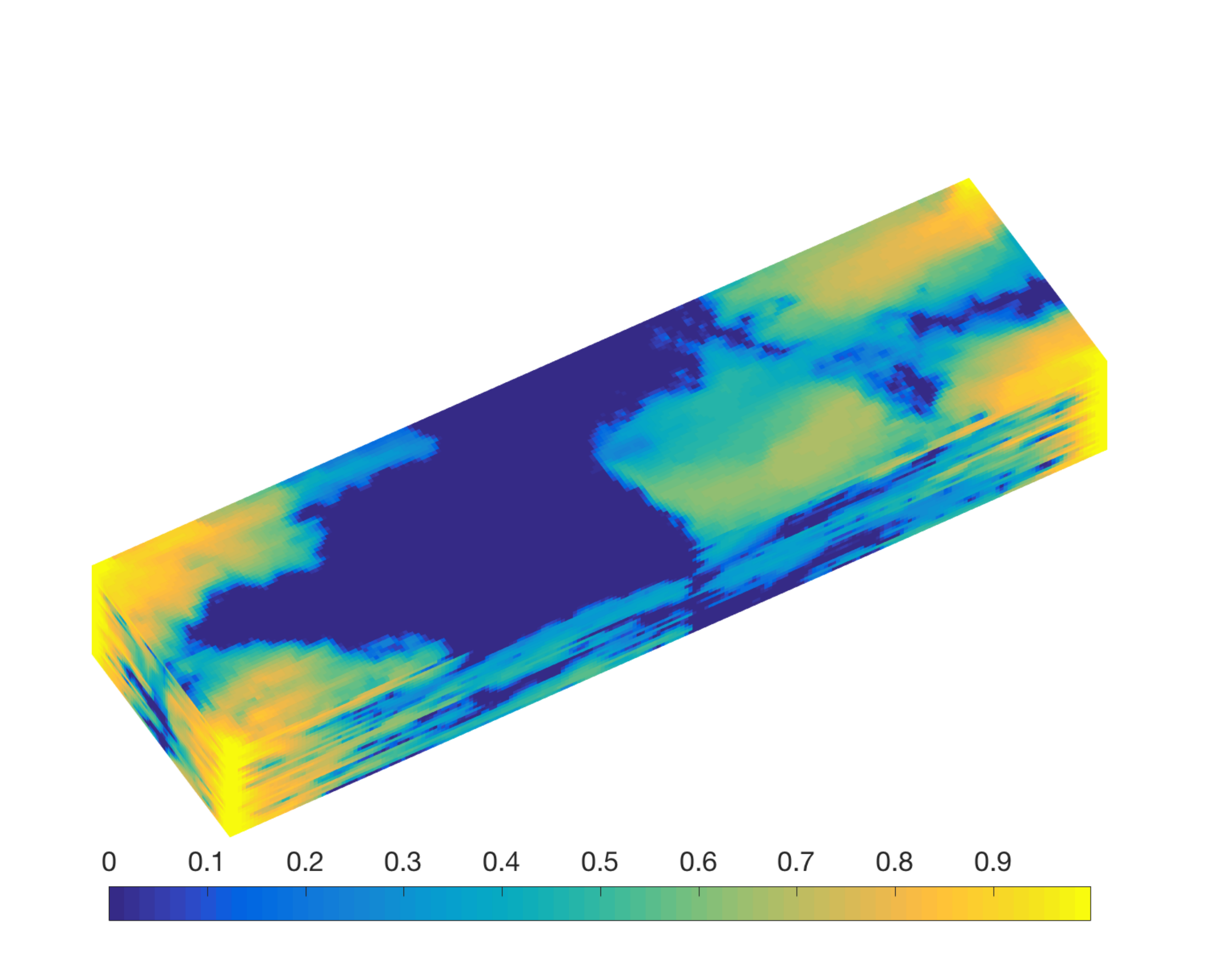}}
	\subfigure[Ms-solution with 1 offline basis and 5 online basis, relative $L^2$ error is 4.8$\%$]{
		\includegraphics[width=3.0in]{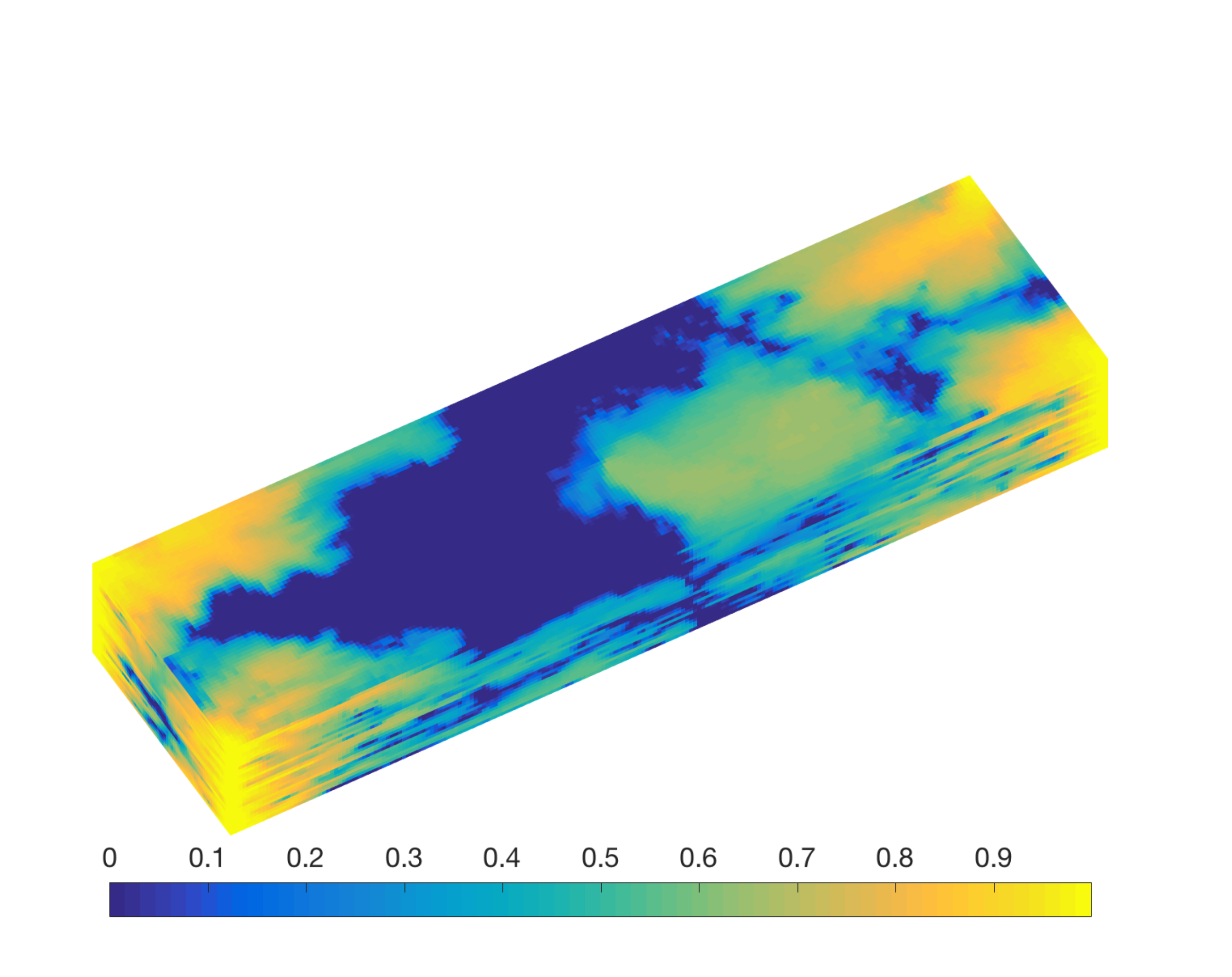}}
		\subfigure[Ms-solution with 1 offline basis and 7 online basis, relative $L^2$ error is 2.3$\%$]{
		\includegraphics[width=3.0in]{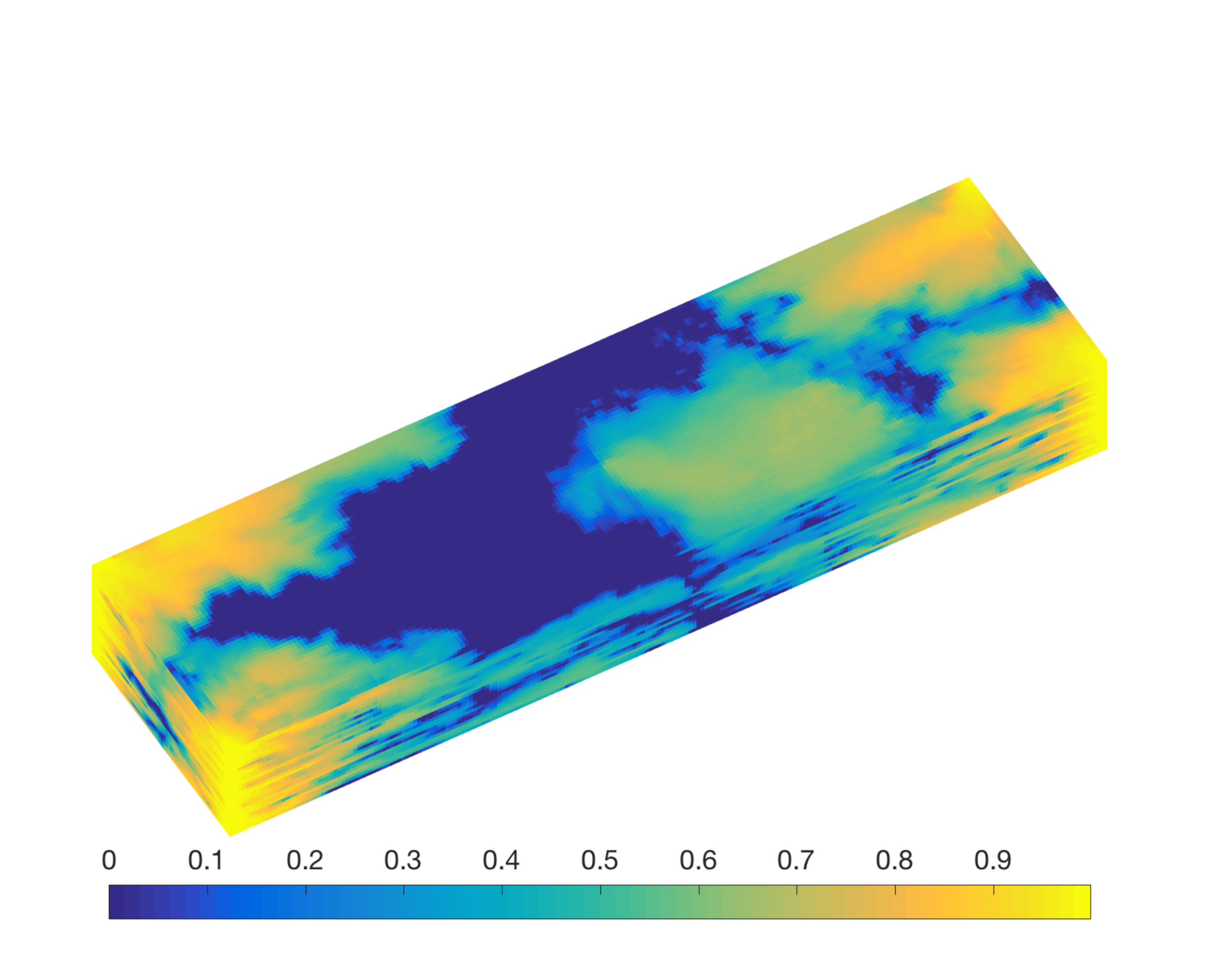}}
	\subfigure[Ms-solution with 4 offline basis and 2 online basis, relative $L^2$ error is 4.3$\%$]{
		\includegraphics[width=3.0in]{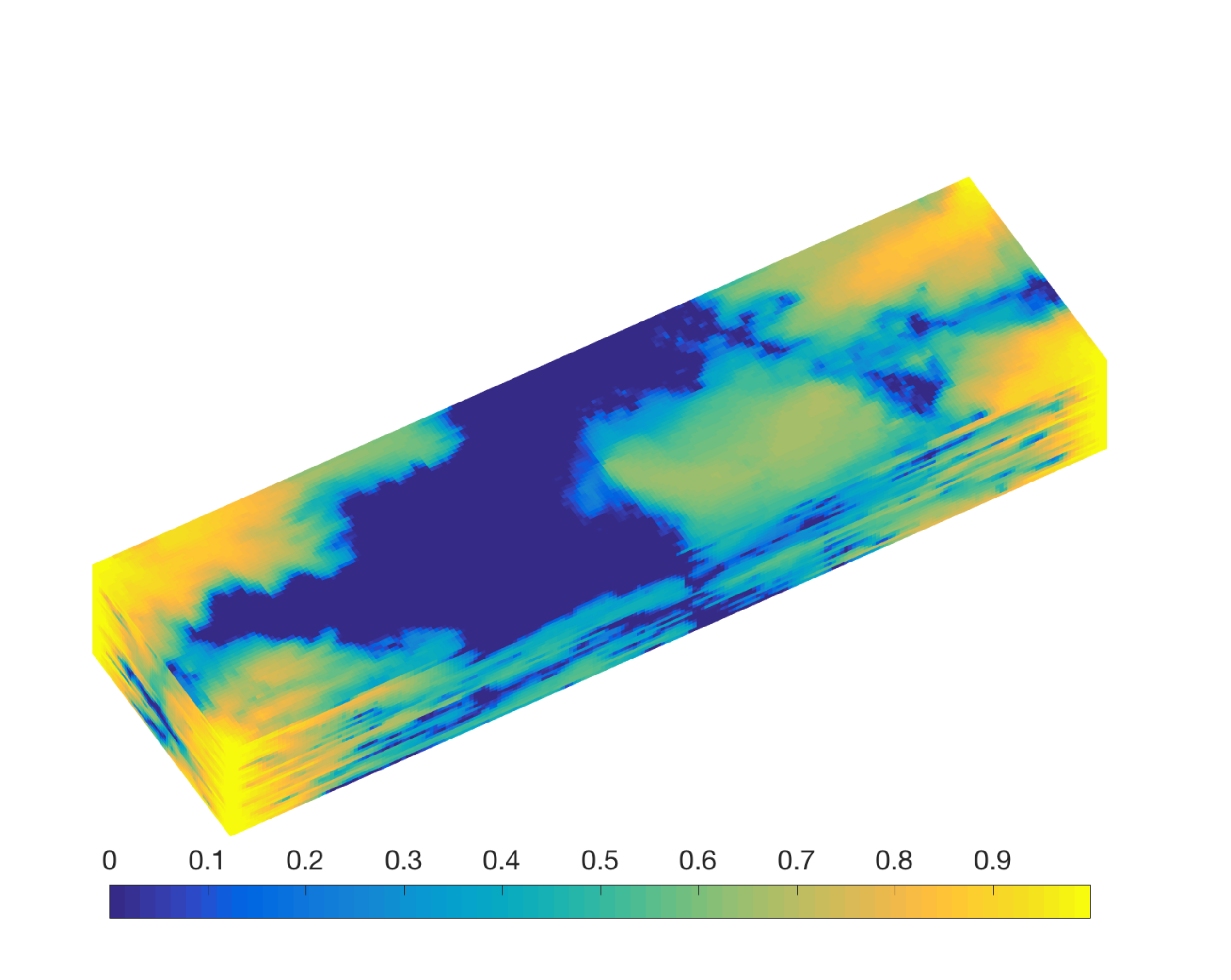}}
		\subfigure[Ms-solution with 4 offline basis and 4 online basis, relative $L^2$ error is 2.3$\%$]{
		\includegraphics[width=3.0in]{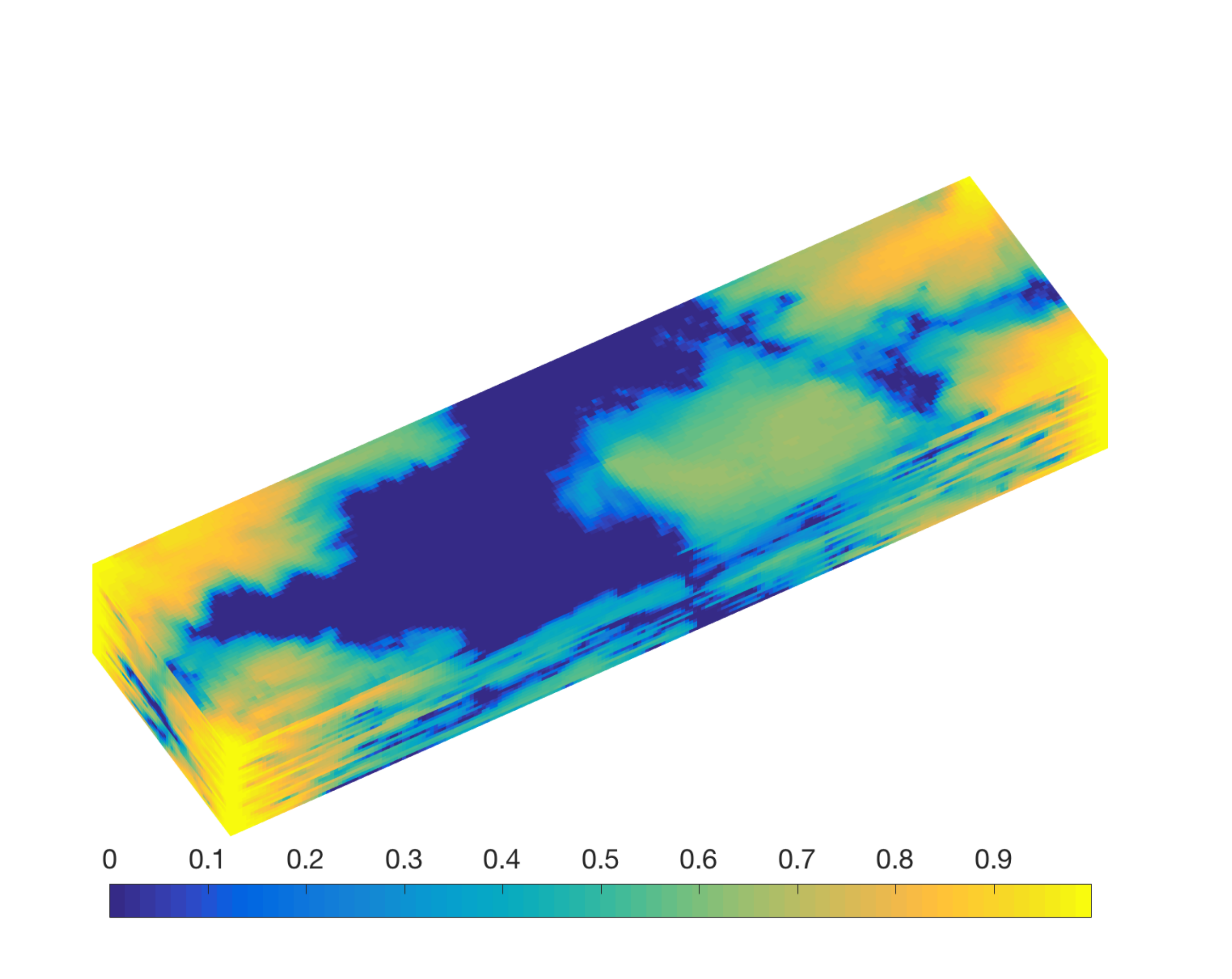}}	
	\caption{Saturation comparison at $t=2500$, for the first 30 layers of SPE 10 as in Figure \ref{model}(b).}
	\label{Sat_model2_5000}
\end{figure}

\begin{figure}[H]
	\centering
	\subfigure[Reference solution]{
		\includegraphics[width=3.0in]{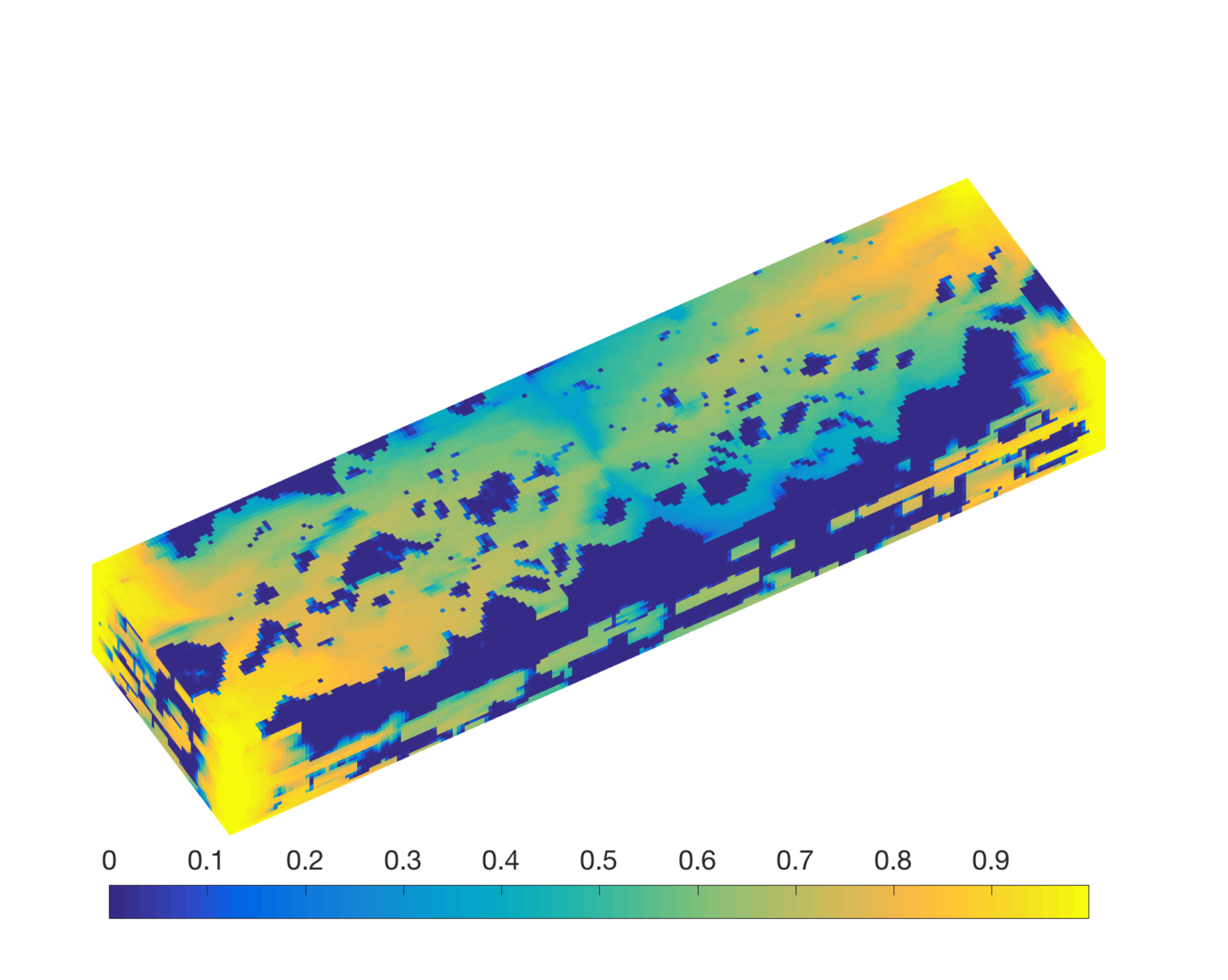}}
	\subfigure[Ms-solution with 1 offline basis and 4 online basis, relative $L^2$ error is 30.7$\%$]{
		\includegraphics[width=3.0in]{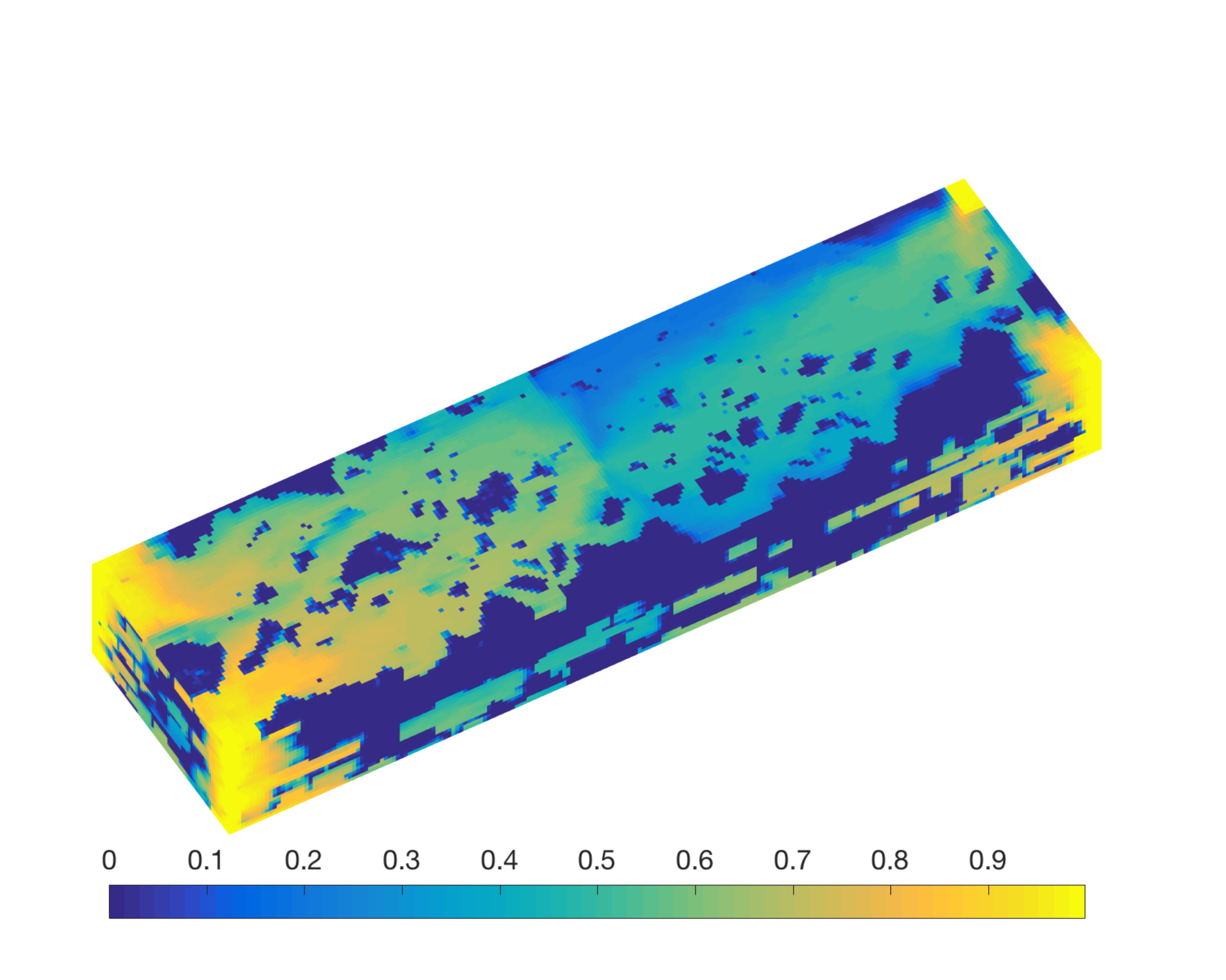}}
		\subfigure[Ms-solution with 1 offline basis and 8 online basis, relative $L^2$ error is 4.9$\%$]{
		\includegraphics[width=3.0in]{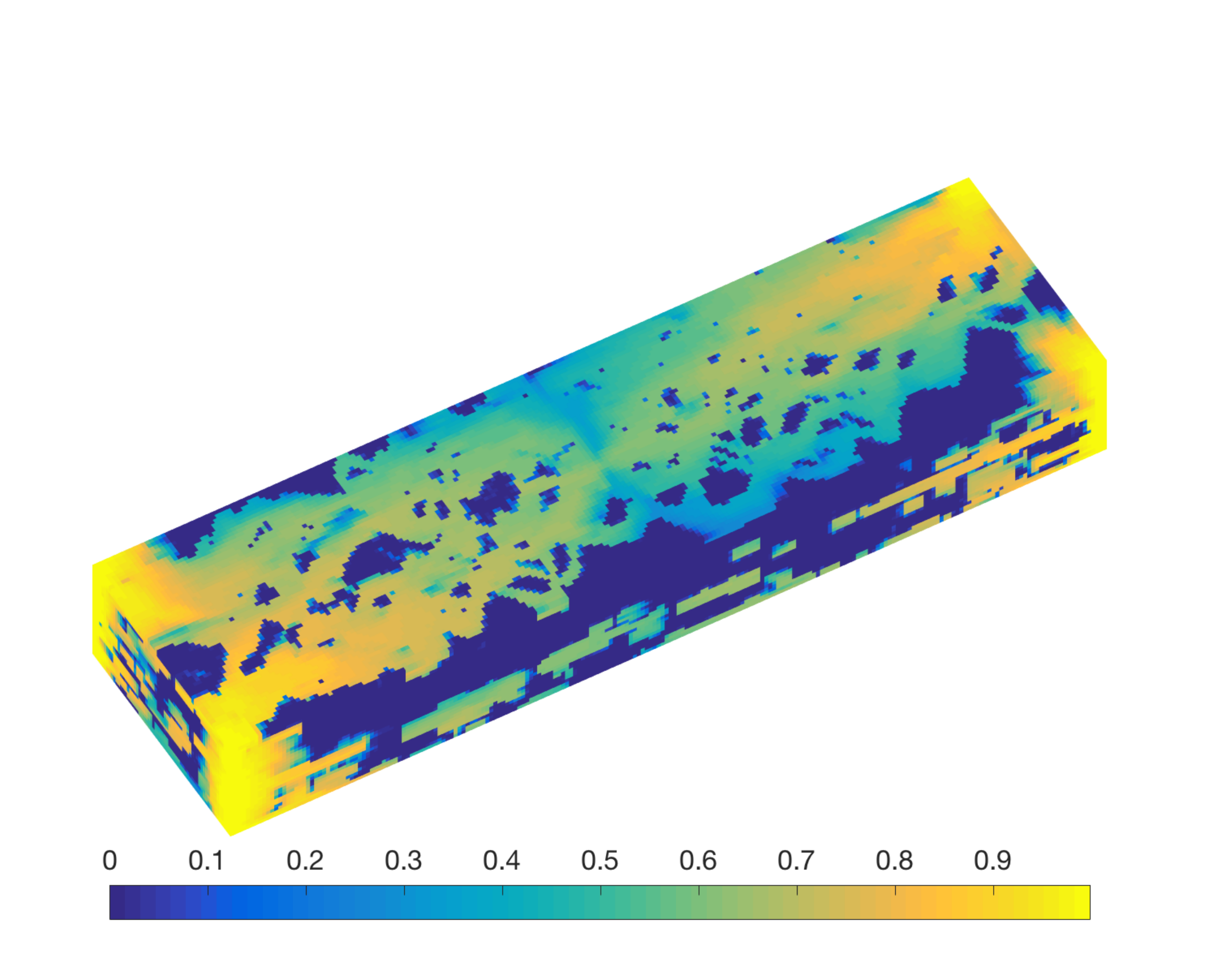}}
	\subfigure[Ms-solution with 1 offline basis and 10 online basis, relative $L^2$ error is 3.4$\%$]{
		\includegraphics[width=3.0in]{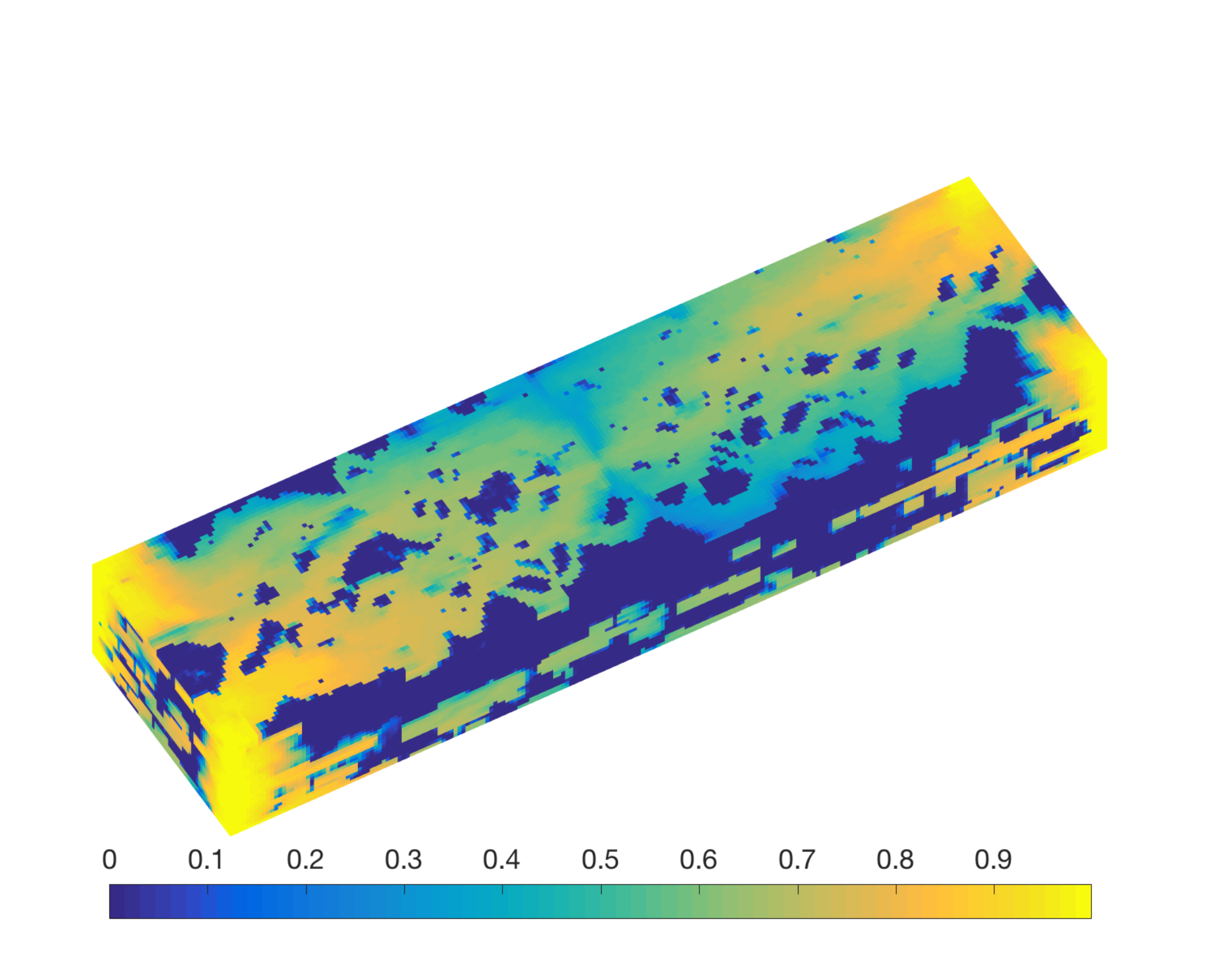}}
	\subfigure[Ms-solution with 4 offline basis and 3 online basis, relative $L^2$ error is 23.2$\%$]{
		\includegraphics[width=3.0in]{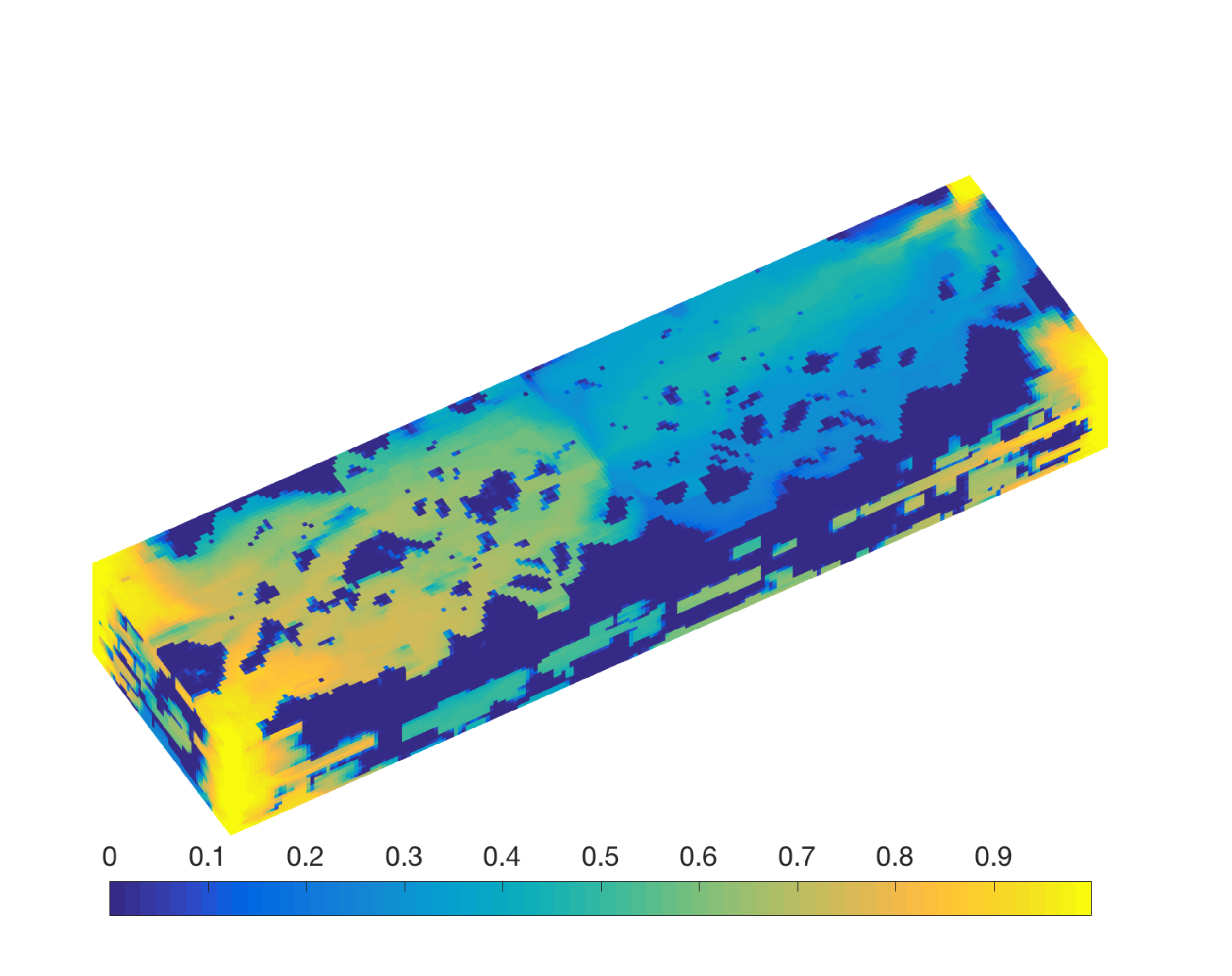}}	
		\subfigure[Ms-solution with 4 offline basis and 7 online basis, relative $L^2$ error is 3.0$\%$]{
		\includegraphics[width=3.0in]{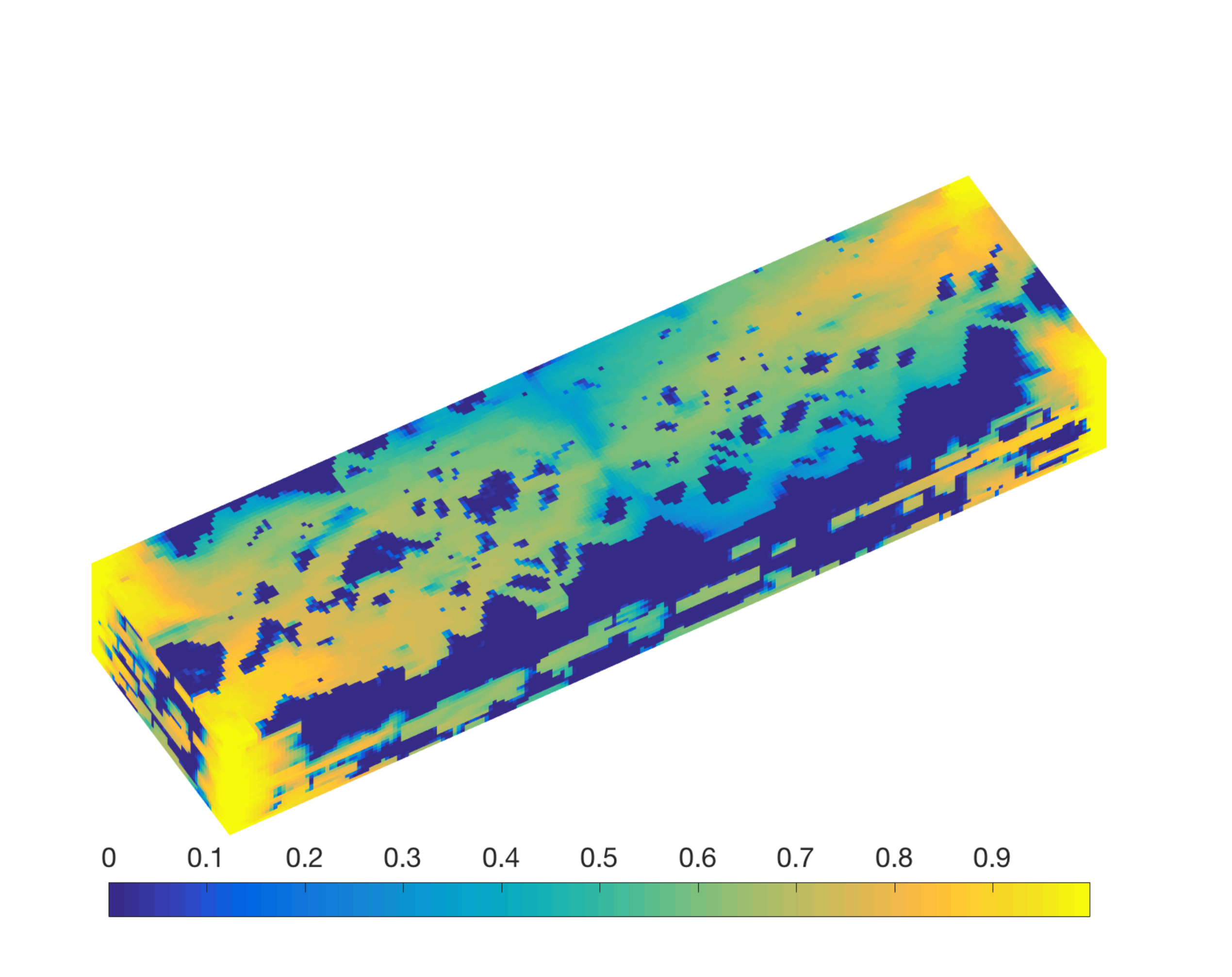}}
	\caption{Saturation comparison at $t=2500$, for the last 30 layers of SPE 10 as in Figure \ref{model}(c).}
	\label{Sat_model3_5000}
\end{figure}

In Figure \ref{Sat_model2_5000}, the saturation plots (at time $t=2500$) for model 2 are given.  Figure \ref{Sat_model2_5000} (a) is the reference solution. Figure \ref{Sat_model2_5000} (b) is the multiscale solution by using 1 offline and 1 online basis functions, which fails to capture much information compared to the reference solution. The relative saturation error is $21.9\%$. After adding 5 online basis functions, the error drops to $4.8\%$, whose saturation plot is given in Figure \ref{Sat_model2_5000} (c).   Figure \ref{Sat_model3_5000} presents the saturation plots (at time $t=2500$) for model 3. Starting with 1 offline basis function, and using 4 online basis functions, the relative saturation error is $30.7\%$, whose corresponding saturation profile is given in Figure \ref{Sat_model3_5000}(b). This saturation profile has large discrepancy compared to the reference one. By adding 8 online basis functions, the error drops to $4.9\%$.  Figure \ref{Sat_model3_5000}(e) is for the case of 4 offline and 3 online basis functions. From 
Figure \ref{Sat_model2_5000} and Figure \ref{Sat_model3_5000}, we see that due to the more complicated feature of the permeability field for model 3, more online basis functions are needed in general to get accurate resutls.

%

%
%

\section{Conclusions}
We have developed an online adaptive multiscale mortar mixed finite element method for flow problems in heterogeneous porous media. We start with a cheap coarse grid solution which are computed by using polynomial functions. The residual from this coarse grid solution is used to compute multiscale basis functions. From the space consisting of both  the new basis functions and the previous one,  we then get a new solution,  then a new residual, and then new basis functions.  We also propose oversampling to compute the online basis functions. We compare the results of oversampling and non-oversampling, which shows that oversampling yields faster convergence speed. An important application of the method is for solving 3D flow and transport problems. We only enrich the solution space of the initial problem, and use this initial solution space for the rest of the simulation along the time. Our numerical results show that the online basis functions produces accurate production file along the time. Our method is efficient and accurate for two-phase flow and transport problem since we do not need to update the multiscale solution space at later time steps.

\section*{Convergence of the algorithm}

In this section, we will give some convergence analysis of the online iterative algorithm. 
We see that the sequence of solutions $\{ \lambda_{\text{ms}}^{l} \}$ generated by our online enrichment algorithm
satisfies a contraction property (equation (\ref{eq:convergence})). Moreover, the convergence rate
is computable, and is related to the residual of the current solution. 
To begin, we prove the following lemma, which gives an a-posteriori error bound for the solution. 
\begin{lemma}
We have
\begin{equation}\label{ineq:posterror}
||\lambda_H^f-\lambda_{\text{ms}}||_{M_H^f}\leq C\displaystyle\sum_{i=1}^{\left\vert{\mathcal{E}_H}\right\vert} ||R_{\omega_i}||_{M_{H,i}^{f,*}}
\end{equation}
where $C$ is a constant that does not depend on the mesh size.
\end{lemma}
Proof: Define $\mathcal{P_H}$ as the $L^2$ projection from $M_H^f$ to $M_H$. Let $\lambda\in M_H^f$ be an arbitrary function in space $M_H^f$. We have
\begin{equation*}
\begin{split}
a_H(\lambda_H^f-\lambda_{\text{ms}},\lambda)=&g_H(\lambda)-a_H(\lambda_{\text{ms}},\lambda)\\
=&g_H(\lambda-\mathcal{P_H} \lambda)+g_H(\mathcal{P_H}\lambda)-a_H(\lambda_{\text{ms}},\lambda)-a_H(\lambda-\lambda_{\text{ms}},\mathcal{P_H}\lambda)\\
=&g_H(\lambda-\mathcal{P_H} \lambda)-a_H(\lambda_{\text{ms}},\lambda-\mathcal{P_H}\lambda)\\
=&\displaystyle\sum_{i=1}^{\left\vert{\mathcal{E}_H}\right\vert}R_{H,i}^f(\lambda-\mathcal{P_H}\lambda)\\
\leq&\displaystyle\sum_{i=1}^{\left\vert{\mathcal{E}_H}\right\vert} ||R_{\omega_i}||_{M_{H,i}^{f,*}}||\lambda-\mathcal{P_H}\lambda||_{M_{H,i}^f}\\
\leq&C\displaystyle\sum_{i=1}^{\left\vert{\mathcal{E}_H}\right\vert} ||R_{\omega_i}||_{M_{H,i}^{f,*}}||\lambda||_{M_H^f}
\end{split}
\end{equation*}
The inequality (\ref {ineq:posterror}) follows by letting $\lambda=\lambda_H^f-\lambda_{\text{ms}}$.

\begin{theorem}
	Using the notation in last subsection, we have
	\begin{equation}\label{eq:convergence}
	||\lambda_H^f-\lambda_{\text{ms}}^{l+1}||_{M_H^f}^2\leq\Bigg(1-\frac{\displaystyle\sum_{i=1}^I||R_{\omega_i}||^2_{M_{H,i}^{f,*}}}{C\displaystyle\sum_{i=1}^{\left\vert{\mathcal{E}_H}\right\vert} ||R_{\omega_i}||_{M_{H,i}^{f,*}}} \Bigg)||\lambda_H^f-\lambda_{\text{ms}}^{l}||_{M_H^f}^2
	\end{equation}
\end{theorem}

Proof: Let $\lambda_H^f$ be the snapshot solution.
From $a_H(\lambda_{\text{ms}}^{l+1},\lambda)=g_H(\lambda)$ for all $\lambda\in M_H^{l+1}$, we can
get $||\lambda_H^f-\lambda_{\text{ms}}^{l+1}||_{M_H^f}^2=\displaystyle{\inf_{\lambda\in M_H^{l+1}}}||\lambda_H^f-\lambda||_{M_H^f}^2$. Let $\lambda_{\text{ms}}^{l+1}=\lambda_{\text{ms}}^{l}+\alpha_1\mu_{H,1}+\alpha_1\mu_{H,2}+\cdots+
\alpha_I\mu_{H,I}$, then
\begin{equation*}
\begin{split}
||\lambda_H^f-\lambda_{\text{ms}}^{l+1}||_{M_H^f}^2\leq&||\lambda_H^f-\lambda_{\text{ms}}^{l}+\alpha_1\mu_{H,1}+\alpha_1\mu_{H,2}+\cdots+
\alpha_I\mu_{H,I}||_{M_H^f}^2\\
=&||\lambda_H^f-\lambda_{\text{ms}}^{l}||_{M_H^f}^2+||\alpha_1\mu_{H,2}+\cdots+
\alpha_I\mu_{H,I}||_{M_H^f}^2\\
&-2a_H(\lambda_H^f-\lambda_{\text{ms}}^{l},\alpha_1\mu_{H,2}+\cdots+
\alpha_I\mu_{H,I})\\
=&||\lambda_H^f-\lambda_{\text{ms}}^{l}||_{M_H^f}^2-\displaystyle\sum_{i=1}^I||R_{\omega_i}||^2_{M_{H,i}^{f,*}}.
\end{split}
\end{equation*}
Then, inequality (\ref{eq:convergence}) can be obtained from the lemma above.

\section*{Acknowledgment}

The research of Eric Chung is partially supported by the Hong Kong RGC General Research Fund (Project 14317516)
and the CUHK Direct Grant for Research 2016-17.

	\bibliographystyle{plain}
	\bibliography{references}

\end{document}